\documentclass[a4paper,reqno]{amsart}
\usepackage[utf8]{inputenc}
\usepackage[english]{babel}
\usepackage{amsmath,amssymb,amsthm,exscale,amsopn,mathtools,color,tikz-cd,fullpage,cite}
\usepackage[colorlinks=true,pdftex,unicode=true,linktocpage,bookmarksopen,hypertexnames=false]{hyperref}
\usepackage{graphicx}

\newtheorem{lem}{Lemma}[section]
\newtheorem{prop}[lem]{Proposition}
\newtheorem{thm}[lem]{Theorem}
\newtheorem{cor}[lem]{Corollary}

\theoremstyle{definition}
\newtheorem{ex}[lem]{Example}
\newtheorem{rem}[lem]{Remark}

\DeclareMathOperator{\Ann}{Ann}
\DeclareMathOperator{\Aut}{Aut}
\DeclareMathOperator{\ch}{char}
\DeclareMathOperator{\clK}{clKdim}

\DeclareMathOperator{\GK}{GKdim}
\DeclareMathOperator{\GL}{GL}
\DeclareMathOperator{\gld}{gldim}
\DeclareMathOperator{\gr}{gr}
\DeclareMathOperator{\h}{ht}
\DeclareMathOperator{\id}{id}

\DeclareMathOperator{\Ker}{Ker}
\DeclareMathOperator{\M}{M}
\DeclareMathOperator{\op}{op}
\DeclareMathOperator{\Qcl}{Q_{cl}}
\DeclareMathOperator{\rk}{rk}
\DeclareMathOperator{\Soc}{Soc}
\DeclareMathOperator{\Span}{Span}
\DeclareMathOperator{\Spec}{Spec}
\DeclareMathOperator{\Sym}{Sym}
\DeclareMathOperator{\Z}{Z}

\newcommand{\mb}{\mathbb}
\newcommand{\mc}{\mathcal}
\newcommand{\mf}{\mathfrak}

\newcommand{\ov}{\overline}
\newcommand{\s}{\subseteq}
\newcommand{\sn}{\subsetneq}
\newcommand{\wh}{\widehat}
\newcommand{\vn}{\varnothing}

\newcommand{\free}[1]{\langle #1 \rangle}

\renewcommand{\le}{\leqslant}
\renewcommand{\ge}{\geqslant}

\title[The structure monoid and algebra \ldots]{The structure monoid and algebra of a non-degenerate set-theoretic solution of the Yang--Baxter equation}

\author[E. Jespers \and {\L}. Kubat \and A. Van Antwerpen]{Eric Jespers \and {\L}ukasz Kubat \and Arne Van Antwerpen}
\address{Department of Mathematics, Vrije Universiteit Brussel, Pleinlaan 2, 1050 Brussel, Belgium}
\email{Eric.Jespers@vub.be}
\email{Lukasz.Kubat@vub.be}
\email{Arne.Van.Antwerpen@vub.be}

\thanks{The first author is supported in part by Onderzoeksraad of Vrije Universiteit Brussel and Fonds voor Wetenschappelijk Onderzoek
(Flanders), grant G016117. The second author is supported by Fonds voor Wetenschappelijk Onderzoek (Flanders), grant G016117.
The third author is supported by Fonds voor Wetenschappelijk Onderzoek (Flanders), via an FWO Aspirant-mandate.}
\subjclass[2010]{Primary: 16N60, 16T25; Secondary: 16R20, 16S36, 16S37}
\keywords{Yang--Baxter equation, set-theoretic solution, prime ideal, (skew) brace}
\date{}

\begin{document}
\begin{abstract}
    For a finite involutive non-degenerate solution $(X,r)$ of the Yang--Baxter equation it is known that the structure monoid
    $M(X,r)$ is a monoid of I-type, and the structure algebra $K[M(X,r)]$ over a field $K$ shares many properties with commutative
    polynomial algebras, in particular, it is a Noetherian PI-domain that has finite Gelfand--Kirillov dimension. In this paper we
    deal with arbitrary finite (left) non-degenerate solutions. Although the structure of both the monoid $M(X,r)$ and the algebra
    $K[M(X,r)]$ is much more complicated than in the involutive case, we provide some deep insights.
    
    In this general context, using a realization of Lebed and Vendramin of $M(X,r)$ as a regular submonoid in the semidirect product 
    $A(X,r)\rtimes\Sym (X)$, where $A(X,r)$ is the structure monoid of the rack solution associated to $(X,r)$, we prove that the
    structure algebra $K[M(X,r)]$ is a module-finite normal extension of a commutative affine subalgebra. In particular, $K[M(X,r)]$
    is a Noetherian PI-algebra of finite Gelfand--Kirillov dimension bounded by $|X|$. We also characterize, in ring-theoretical terms of
    $K[M(X,r)]$, when $(X,r)$ is an involutive solution. This characterization provides, in particular, a positive answer to the
    Gateva-Ivanova conjecture concerning cancellativity of $M(X,r)$. 
    
    These results allow us to control the prime spectrum of the algebra $K[M(X,r)]$ and to describe the Jacobson radical and prime
    radical of $K[M(X,r)]$. Finally, we give a matrix-type representation of the algebra $K[M(X,r)]/P$ for each prime ideal
    $P$ of $K[M(X,r)]$. As a consequence, we show that if $K[M(X,r)]$ is semiprime then there exist finitely many finitely generated
    abelian-by-finite groups, $G_1,\dotsc,G_m$, each being the group of quotients of a cancellative subsemigroup of $M(X,r)$ such that
    the algebra $K[M(X,r)]$ embeds into $\M_{v_1}(K[G_1])\times\dotsb\times\M_{v_m}(K[G_m])$, a direct product of matrix algebras.
\end{abstract}
\maketitle

\section*{Introduction}

Let $V$ be a vector space over a field $K$. A linear map $R\colon V\otimes V\to V\otimes V$ is called a solution of the Yang--Baxter equation
(or braided equation) if \[(R\otimes{\id})\circ({\id}\otimes R)\circ(R\otimes{\id})=({\id}\otimes R)\circ(R\otimes{\id})\circ({\id}\otimes R).\] 
Recall that this  equation originates from papers by Baxter \cite{Ba} and Yang \cite{Ya} on statistical physics and the search for solutions has
attracted numerous studies both in mathematical physics and pure mathematics.

As the study of arbitrary solutions is complex, Drinfeld, in 1992 \cite{Dr}, proposed to study the solutions which are induced
by a linear extension of a map $r\colon X\times X\to X\times X$, where $X$ is a basis of $V$. In this case $r$ satisfies
\[(r\times{\id})\circ({\id}\times r)\circ(r\times{\id})=({\id}\times r)\circ(r\times{\id})\circ({\id}\times r),\]
and one says that $(X,r)$ is a set-theoretic solution of the Yang--Baxter equation.
For any $x,y \in X$, we put $r(x,y)=(\lambda_x(y),\rho_y(x))$. Since that late 1990's several ground-breaking results were discovered
on this topic, including these by Gateva-Ivanova and Van den Bergh \cite{GVdB}, Etingof, Schedler, and Soloviev \cite{ESS} and
Lu, Yan, and Zhu \cite{LYZ}. The investigations on the subject have intensified even more since the discovery
of several algebraic structures associated to set-theoretic solutions. A particular nice class of set-theoretic solutions
$(X,r)$ are the bijective (i.e., $r$ is a bijection) solutions that are left and right non-degenerate (i.e., each $\lambda_x$,
respectively each $\rho_x$, is a bijection). If furthermore $r^2=\id$ then the solution is said to be involutive. In order to
deal with such involutive solutions Rump \cite{R1,R2} introduced the new algebraic structure called ``(left) brace'' and Guarnieri
and Vendramin \cite{GV} extended this algebraic structure to a ``(left) skew brace'' in order to also deal with arbitrary bijective
non-degenerate solutions. Many fundamental results on these structures already have have been obtained 
\cite{Bach2,Bach,BCJO,CGS,CJO,Ga1,KSV,LV2,R3,SV}. In particular, it has been shown that determining all finite
(i.e., $X$ is a finite set) bijective involutive non-degenerate solutions is equivalent
to describing all finite (left) braces. In \cite{BCJ} a concrete realization of this description has been given. Moreover, braces have lent 
themselves as a novel method to solve questions in group and ring theory. For instance, Amberg, Dickenschied, and Sysak in \cite{ADS} posed the 
question whether the adjoint group of a nil ring is an Engel group, and Zelmanov asked a similar question in the context of nil algebras over an 
uncountable field. Smoktunowicz, using tools related to braces, gave negative answers to both of these questions in \cite{Sm}. Also non-bijective
set-theoretic solutions are of importance and receive attention. For example Lebed in \cite{Le} shows that idempotent solutions
provide a unified treatment of factorizable monoids, free and free commutative monoids, distributive lattices and Young tableaux
and Catino, Colazzo, and Stefanelli \cite{CCS}, and Jespers and Van Antwerpen \cite{JVA} introduced the algebraic
structure called ``(left) semi-brace'' to deal with solutions that are not necessarily non-degenerate or that are idempotent.

In \cite{ESS} Etingof, Schedler, and Soloviev and Gateva-Ivanova and Van den Bergh in \cite{GVdB} introduced the following associated algebraic 
structures to a set-theoretic solution $(X,r)$: the structure group $G(X,r)=\gr(X\mid xy=uv\text{ if }r(x,y)=(u,v))$ and the structure monoid
$M(X,r)=\free{X \mid xy=uv\text{ if }r(x,y)=(u,v)}$. In case $(X,r)$ is finite involutive and non-degenerate it is shown that the group $G(X,r)$ is 
solvable and it is naturally embedded into the semidirect product $\mb{Z}^{(X)}\rtimes\Sym(X)$, where $\Sym (X)$ acts naturally on the free
abelian group $\mb{Z}^{(X)}$ of rank $|X|$. It turns out that $G(X,r)=\gr((x,\lambda_x)\mid x\in X)$ and, in particular, these groups are (free 
abelian)-by-finite. Furthermore, in \cite{JO2} it is shown that $M(X,r)$ is embedded in $G(X,r)$ and the latter is the group of fractions of 
$M(X,r)$. Furthermore, $G(X,r)$ and $\mc{G}(X,r)=\gr(\lambda_x\mid x\in X)$ are left braces and, for finite $X$, groups of the type $\mc{G}(X,r)$ 
correspond to all finite (left) braces (for details we refer to \cite{Ce}).
In \cite{GVdB} Gateva-Ivanova and Van den Bergh showed that these structure groups are groups of I-type and, in particular, they are finitely 
generated and torsion-free, i.e., Bieberbach groups. These groups and monoids are of combinatorial interest and their associated monoid algebra 
$K[M(X,r)]$, simply called the structure algebra of $(X,r)$ (as it is the algebra generated by the set $X$ and with defining relations
$xy=uv\text{ if }r(x,y)=(u,v)$), provide non-trivial examples of quadratic algebras. That is, they are positively graded algebras generated by the 
homogeneous part of degree $1$ and with defining homogeneous degree $2$ relations. The structure algebras have similar homological properties to 
polynomial algebras in finitely many commuting variables, in particular they are Noetherian domains that satisfy a polynomial identity (PI-algebras) 
and have finite Gelfand--Kirillov dimension. 

Recently, Lebed and Vendramin \cite{LV} studied the structure group $G(X,r)$ for arbitrary finite bijective non-degenerate solutions
(i.e., not necessarily involutive). In \cite{LV,LYZ,So} they associate, via a bijective $1$-cocycle, to the structure group $G(X,r)$
the structure group $G(X,\triangleleft_r)$ of the structure rack $(X,\triangleleft_r)$ of $(X,r)$. 
As a consequence, it follows that again the groups $G(X,r)$ are 
abelian-by-finite. Recall that a set $X$ with a self-distributive operation $\triangleleft$ is called a rack if the map $y\mapsto y\triangleleft x$
is bijective, for any $x\in X$ (cf. \cite{Ka}). In contrast to the involutive case, the set $X$ is not necessarily canonically embedded into
$G(X,r)$, the reason being that $M(X,r)$ need not be cancellative in general (i.e., it is not necessarily embedded in a group). Hence, for
an arbitrary solution $(X,r)$ the structure monoid $M(X,r)$ contains more information on the original solution. However, it is in general not
true that two set-theoretic solutions $(X,r)$ and $(Y,s)$ are isomorphic if and only if the monoids $M(X,r)$ and $M(Y,s)$ are isomorphic.
This does hold if one of both solutions (and thus both) is assumed to be an involutive non-degenerate set-theoretic solution.

In this paper we give a structural approach of the study of the structure monoid $M(X,r)$ and the structure algebra $K[M(X,r)]$ for an arbitrary 
bijective left non-degenerate solution $(X,r)$. In the same spirit as in \cite{LV}, in Section~\ref{sec:1} we associate a structure monoid, called
the derived structure monoid and denoted $A(X,r)$, to such a solution and we show that the monoid $M(X,r)$ is a regular submonoid of
$A(X,r)\rtimes\Sym(X)$, i.e., there is a bijective $1$-cocycle $M(X,r)\to A(X,r)$. Again $A(X,r)$ turns out to be the structure monoid of a rack. This 
description allows us to study, in Sections \ref{sec:2} and \ref{sec:4}, the algebraic structure of the monoids $A(X,r)$ and $M(X,r)$ and the
structure algebras $K[A(X,r)]$ and $K[M(X,r)]$. It is shown that for a finite bijective left non-degenerate solution $(X,r)$ the monoid $A(X,r)$ 
(respectively $M(X,r)$) is central-by-finite (respectively abelian-by-finite), i.e., they are finite ``modules'' over finitely generated central
(respectively commutative) submonoids. Hence, the both structure algebras are Noetherian and PI. Furthermore, these algebras are closely related to
polynomial algebras in finitely many commuting variables, for instance we show that the classical Krull dimensions and Gelfand--Kirillov dimensions
of both algebras $K[A(X,r)]$ and $K[M(X,r)]$ coincide and are equal to $\rk A(X,r)=\rk M(X,r)$, i.e., the rank of the respective monoids (that is the 
largest possible rank of a free abelian submonoid). Moreover, this dimension is bounded by $|X|$ and it also is shown that these dimensions are 
determined by the orbits of subsolutions of the rack solution $(X,s)$ associated to $(X,r)$. Gateva-Ivanova in \cite{Ga2} conjectured that the 
structure monoid of a finite square-free (i.e., $r(x,x)=(x,x)$ for all $x\in X$) non-degenerate solution $(X,r)$ is cancellative if and only if the 
solution $(X,r)$ is involutive. Using the structural results we prove that this conjecture is true, even without assuming that the solution $(X,r)$
is square-free. Moreover, we show that the involutiveness of a solution is characterized by many properties of the structure algebra $K[M(X,r)]$.
Among others, this coincides with the maximality of the Gelfand--Kirillov dimension, i.e., $\GK K[M(X,r)]=|X|$, and it is equivalent with $K[M(X,r)]$ 
being a prime algebra or a domain.

In Section \ref{sec:3} we study the prime ideals of the monoid $A(X,r)$ and the prime ideals of the related algebra. It is shown that prime
ideals of $A(X,r)$ are in correspondence with specific subsolutions of the rack solution $(X,s)$ associated to $(X,r)$.
Furthermore, we provide a description of the prime ideals of $K[A(X,r)]$. 

In Section \ref{sec:5} we study the prime ideals of the monoid $M(X,r)$ and the prime ideals of $K[M(X,r)]$. In \cite{GJ} prime ideals
of monoids of IG-type were studied by Goffa and Jespers. It is shown that for a finite left non-degenerate solution $(X,r)$
the prime ideals of $A(X,r)$ determine the prime ideals of $M(X,r)$; these results are similar to those obtained for monoids of IG-type,
i.e., regular submonoids of the holomorph of a finitely generated cancellative abelian monoid. Furthermore, prime ideals of the algebra
$K[M(X,r)]$, where $(X,r)$ is a finite bijective left non-degenerate solution, which intersect the monoid trivially correspond to prime
ideals of the group algebra $K[G(X,r)]$. As $G(X,r)$ is a finitely generated finite-conjugacy group (FC-group for short) the prime ideals
of $K[G(X,r)]$ are easy to describe. For more fundamental results of prime ideals of finitely generated abelian-by-finite groups, or more
general, polycyclic-by-finite groups, we refer the reader to the fundamental work of Roseblade \cite{Ro}. 

In \cite{GJO} Gateva-Ivanova, Jespers, and Okni\'nski and, in \cite{JOVC,JVC}, Jespers, Okni\'nski, and Van Campenhout studied
the prime ideals of quadratic algebras coming from monoids of quadratic type, these are monoids defined on a finite set $X$ of cardinality
$n$ and defined by $\binom{n}{2}$ monomial relations of degree two so that the associated map $r\colon X\times X\to X\times X$ is non-degenerate;
but it does not have to be a set-theoretic solution of the Yang--Baxter equation. They showed that the intersection of such prime ideals
with the monoid is highly dependent on the divisibility structure of the monoid. In Section~\ref{sec:6} the divisibility structure of $M(X,r)$
is studied. It is shown that the intersection of a prime ideal of $K[M(X,r)]$ with $M(X,r)$ is determined by divisibility properties. These results 
allow to give a description of the Jacobson radical $\mc{J}(K[M(X,r)])$ and prime radical $\mc{B}(K[M(X,r)])$ of $K[M(X,r)]$.

In the final Section~\ref{sec:7} we prove a matrix-type representation of the prime algebra $K[M(X,r)]/P$ for each prime ideal $P$ of
$K[M(X,r)]$. It is shown that the classical ring of quotients $\Qcl(K[M(X,r)]/P)$ of $K[M(X,r)]/P$ is the same as $\Qcl(\M_v(K[G]/P_0))$,
where $P_0$ is a prime ideal of a group algebra $K[G]$ with $G$ the group of quotients of a cancellative subsemigroup of $M(X,r)$ and
$v\ge 1$ is determined by the number of orthogonal cancellative subsemigroups of an ideal in $M(X,r)/(P\cap M(X,r))$. As a consequence,
we show that if, furthermore, $K[M(X,r)]$ is semiprime then there exist finitely many finitely generated abelian-by-finite groups, 
say $G_1,\dotsc,G_m$, each being the group of quotients of a cancellative subsemigroup of $M(X,r)$, such that
$K[M(X,r)]$ embeds into $\M_{v_1}(K[G_1])\times\dotsb\times\M_{v_m}(K[G_m])$ for some $v_1,\dotsc,v_m\ge 1$.

\section{Preliminaries}\label{sec:1}
Let $X$ be a non-empty set and $r\colon X\times X\to X\times X$ a map denoted as
\[r(x,y)=(\lambda_x(y),\rho_y(x))\] for $x,y\in X$. Then $(X,r)$ is a solution of the Yang--Baxter equation if and only if,
for any $x,y,z\in X$, the following equalities hold:
\begin{align}
	\lambda_x(\lambda_y(z)) &=\lambda_{\lambda_x(y)}(\lambda_{\rho_y(x)}(z)),\label{YB1}\\
	\lambda_{\rho_{\lambda_y(z)}(x)}(\rho_z(y)) &=\rho_{\lambda_{\rho_y(x)}(z)}(\lambda_x(y)),\label{YB2}\\
	\rho_z(\rho_y(x)) &=\rho_{\rho_z(y)}(\rho_{\lambda_y(z)}(x)).\nonumber
\end{align}
For a solution $(X,r)$ we define its structure monoid (we use the terminology introduced in
\cite{ESS}; in \cite{Ga2} this is called the monoid associated with $(X,r)$)
\[M(X,r)=\free{X\mid xy=\lambda_x(y)\rho_y(x)\text{ for all }x,y\in X}.\] It turns out that in the study of $M(X,r)$ the
derived structure monoid (we use terminology similar as in \cite{So} in the context of groups)
\[A(X,r)=\free{X\mid x\lambda_x(y)=\lambda_x(y)\lambda_{\lambda_x(y)}(\rho_y(x))\text{ for all }x,y\in X}\]
plays a crucial role. If we put $z=\lambda_x(y)$ then the defining relations of $A(X,r)$ can be rewritten as $xz=z\sigma_z(x)$, where
$\sigma_z(x)=\lambda_z(\rho_{\lambda_x^{-1}(z)}(x))$. Hence, \[A(X,r)=\free{X\mid xz=z\sigma_z(x)\text{ for all }x,z\in X}.\]
Moreover, if $(X,r)$ is bijective left non-degenerate, it can be proved that $(X,r^{-1})$ is automatically a left non-degenerate solution.
In this case, writing $r^{-1}(x,y)=(\hat{\lambda}_x(y),\hat{\rho}_y(x))$ for $x,y\in X$, it can be verified that
\begin{equation}\label{sig}
	\sigma_z(x)=\lambda_z(\rho_{\lambda_x^{-1}(z)}(x))=\lambda_z(\hat{\lambda}_z^{-1}(x))
\end{equation}
for all $x,z\in X$. Notice that the second equality in \eqref{sig} leads to $\sigma_z\in\Sym(X)$. Note also that if the solution $(X,r)$ is involutive,
then $\sigma_x=\id$ for all $x\in X$ and thus $A(X,r)$ is the free abelian monoid of rank $|X|$.

Since the defining relations of $M(X,r)$ and $A(X,r)$ are homogeneous, both these monoids inherit a gradation determined by the length
function on words in the free monoid on $X$. We shall freely use this fact throughout the paper. Moreover, the length of an element $s$
in one of the monoids under consideration will be denoted by $|s|$.

Let $(X,r)$ and $(Y,s)$ be solutions of the Yang--Baxter equation. We say that a map $f\colon X\to Y$
is a morphism of solutions (and we write $f\colon(X,r)\to(Y,s)$) if $(f\times f)\circ r=s\circ(f \times f)$
or, in other words, if the following diagram
\[\begin{tikzcd}[row sep=30pt,column sep=30pt]
    X\times X\ar[r,"f\mathop{\times}f"]\ar[d,"r",swap] & Y\times Y\ar[d,"s"]\\
    X\times X\ar[r,"f\mathop{\times}f"] & Y\times Y
\end{tikzcd}\]
is commutative. Moreover, the solutions $(X,r)$ and $(Y,s)$ are called isomorphic provided there exists
a bijective morphism of solutions $f\colon(X,r)\to(Y,s)$. Two involutive non-degenerate solutions $(X,r)$ and $(Y,s)$
are isomorphic if and only if their structure monoids $M(X,r)$ and $M(Y,s)$ are isomorphic. To see this, it is sufficient
to induce an action of $r$ on the words of length two in the alphabet $X$ and observe that these orbits are of size two or smaller.
However, the following example shows that this is no longer true for non-involutive solutions.

\begin{ex}\label{ex:1}
    Let $X=\{x_1,x_2,x_3\}$. Define $\sigma_1=(2,3)$, $\sigma_2=(1,3)$, $\sigma_3=(1,2)$ and consider the maps
    $r,s\colon X\times X\to X\times X$ given by
    \[r(x_i,x_j)=(x_j,x_{\sigma_j(i)})\qquad\text{and}\qquad s(x_i,x_j)=(x_{\sigma_i(j)},x_i).\]
    It is easy to check that both $(X,r)$ and $(X,s)$ are bijective (in fact, $r^3=s^3=\id$) non-degenerate
    solutions of the Yang--Baxter equation. Moreover, $M(X,r)=A(X,r)=A(X,s)=M(X,s)$. However, $(X,r)$ and $(X,s)$
    are not isomorphic as solutions. Indeed, if $f\colon(X,r)\to(X,s)$ were an isomorphism of solutions then,
    in particular, $f\circ\sigma_x=f$ for all $x\in X$, which would lead to $\sigma_x=\id$, a contradiction.
\end{ex}

The remaining part of this section is based on the work of Lebed and Vendramin \cite{LV}. For completeness' sake and to translate  their results
on bijective $1$-cocycles into the language of regular submonoids, which will be crucial to all sections in this paper,  we include detailed proofs.  

By an action of a monoid $M$ on a monoid $A$ we mean a left action by automorphisms, that is a morphism of monoids
$\theta\colon M\to\Aut(A)$ (multiplication in $\Aut(A)$ will be often written as a juxtaposition). Recall that a map
$\varphi\colon M\to A$ is called a bijective $1$-cocycle with respect to the action $\theta$ provided $\varphi$ is bijective,
$\varphi(1)=1$ (i.e., $\varphi$ preserves units of monoids) and satisfies the $1$-cocycle condition 
\[\varphi(xy)=\varphi(x)\theta(x)(\varphi(y))\] for all $x,y\in M$.

\begin{lem}\label{lem:1}
	Assume that $\theta\colon M\to\Aut(A)$ is an action and $\varphi\colon M\to A$ is a bijective $1$-cocycle with respect to $\theta$.
	For a congruence $\eta$ on $M$ define \[\varphi(\eta)=\{(\varphi(x),\varphi(y)):(x,y)\in\eta\}\s A\times A.\]
	If the congruence $\eta$ satisfies the following properties
	\begin{enumerate}
		\item $\eta\s\Ker\theta=\{(x,y)\in M\times M:\theta(x)=\theta(y)\}$ and 
		\item $\varphi(\eta)=\{(\theta(z)(\varphi(x)),\theta(z)(\varphi(y))):(x,y)\in\eta\}$ for all $z\in M$
	\end{enumerate}
	then $\varphi(\eta)$ is a congruence on $A$. Moreover, $\theta$ induces an action $\ov{\theta}\colon M/\eta\to\Aut(A/\varphi(\eta))$
	and $\varphi$ induces a bijective $1$-cocycle $\ov{\varphi}\colon M/\eta\to A/\varphi(\eta)$ with respect to $\ov{\theta}$.
	\begin{proof}
		Using bijectivity of $\varphi$ it is easy to verify that $\varphi(\eta)$ is an equivalence relation on $A$.
		To check that $\varphi(\eta)$ is a left congruence fix $(a,b)\in\varphi(\eta)$ and $c\in A$. Since $\varphi$ is bijective,
		we can write $c=\varphi(z)$ for some $z\in M$. By (2) we get $a=\theta(z)(\varphi(x))$ and $b=\theta(z)(\varphi(y))$ for some
		$(x,y)\in\eta$. Now \[ca=\varphi(z)\theta(z)(\varphi(x))=\varphi(zx)\qquad\text{and}\qquad
		cb=\varphi(z)\theta(z)(\varphi(y))=\varphi(zy).\]
		Because $\eta$ is a left congruence we get $(zx,zy)\in\eta$, and thus $(ca,cb)\in\varphi(\eta)$.
		To prove that $\varphi(\eta)$ is a right congruence assume that $(a,b)\in\varphi(\eta)$ and $c\in A$. 
		By the definition of $\varphi(\eta)$ there exists $(x,y)\in\eta$ such that $(a,b)=(\varphi(x),\varphi(y))$. By (1) we know that
		$\theta(x)=\theta(y)$. Moreover, bijectivity of $\varphi$ assures that $c=\theta(x)(\varphi(z))=\theta(y)(\varphi(z))$
		for some $z\in M$. Now \[ac=\varphi(x)\theta(x)(\varphi(z))=\varphi(xz)\qquad\text{and}\qquad
		bc=\varphi(y)\theta(y)(\varphi(z))=\varphi(yz).\]
		Since $\eta$ is a right congruence we get $(xz,yz)\in\eta$. Hence $(ac,bc)\in\varphi(\eta)$.
		
		To finish the proof observe that (1) implies that there exists an action of $M/\eta$ on $A$ induced by $\theta$.
		Moreover, (2) guarantees that the latter action induces an action $\ov{\theta}\colon M/\eta\to\Aut(A/\varphi(\eta))$.
		Finally, it is clear that $\varphi$ induces a map $\ov{\varphi}\colon M/\eta\to A/\varphi(\eta)$ satisfying
		$\ov{\varphi}(1)=1$ and the cocycle condition with respect to $\ov{\theta}$. Moreover,
		bijectivity of $\ov{\varphi}$ follows easily from bijectivity of $\varphi$.
	\end{proof}
\end{lem}

\begin{lem}\label{lem:2}
	Assume that $\theta\colon M\to\Aut(A)$ is an action and $\varphi\colon M\to A$ is a bijective $1$-cocycle with respect to $\theta$.
	Let $\mc{G}=\theta(M)\s\Aut(A)$, which is a submonoid of $\Aut(A)$. Then the map $f\colon M\to A\rtimes\mc{G}$
	defined as $f(x)=(\varphi(x),\theta(x))$ for $x\in M$
	is an injective morphism of monoids. In particular, \[M\cong f(M)=\{(a,\phi(a)):a\in A\}\s A\rtimes\mc{G},\] where the map 
	$\phi=\theta\circ\varphi^{-1}\colon A\to\mc{G}$ satisfies $\phi(a)\phi(b)=\phi(a\phi(a)(b))$ for $a,b\in A$.
	\begin{proof}
		We have $f(1)=(\varphi(1),\theta(1))=(1,\id)$. Moreover, if $x,y\in M$ then
		\begin{align*}
			f(xy) &=(\varphi(xy),\theta(xy))\\
				  &=(\varphi(x)\theta(x)(\varphi(y)),\theta(x)\theta(y))\\
				  &=(\varphi(x),\theta(x))(\varphi(y),\theta(y))\\
				  &=f(x)f(y).
		\end{align*}
		Since $\varphi$ is injective, $f$ is injective as well. Finally, if $a,b\in A$ then $a=\varphi(x)$ and $b=\varphi(y)$ for some $x,y\in M$.
		Therefore \[\phi(a)\phi(b)=\theta(x)\theta(y)=\theta(xy)=\phi(\varphi(xy))=\phi(\varphi(x)\theta(x)(\varphi(y)))=\phi(a\phi(a)(b)).\]
		Hence the result follows.
	\end{proof}
\end{lem}

\begin{prop}\label{prop:1}
	Assume that $(X,r)$ is a left non-degenerate solution of the Yang--Baxter equation. Let $A=A(X,r)$, $M=M(X,r)$
	and $\mc{G}=\mc{G}(X,r)=\gr(\lambda_x\mid x\in X)\s\Sym(X)$.
	\begin{enumerate}
		\item There exists an action $\theta\colon M\to\Aut(A)$ and a bijective $1$-cocycle $\varphi\colon M\to A$
		with respect to $\theta$ satisfying $\theta(x)=\lambda_x$ and $\varphi(x)=x$ for $x\in X$. In particular, $\mc{G}=\theta(M)$.
		\item The map $f\colon M\to A\rtimes\mc{G}$ defined as $f(x)=(\varphi(x),\theta(x))$ for $x\in M$ is an
		injective morphism of monoids. In particular, \[M\cong f(M)=\{(a,\phi(a)):a\in A\}\s A\rtimes\mc{G},\]
		where the map $\phi=\theta\circ\varphi^{-1}\colon A\to\mc{G}$ satisfies $\phi(a)\phi(b)=\phi(a\phi(a)(b))$ for $a,b\in A$.
		That is, $M$ is a regular submonoid of the semidirect product $A\rtimes\mc{G}$.
		\item If the set $X$ is finite then $\mc{G}$ is a finite group.
	\end{enumerate}
	\begin{proof}
		Let $F$ denote the free monoid on $X$. Define the action $\vartheta\colon F\to\Aut(F)$ by the rule $\vartheta(x)=\lambda_x$
		for $x\in X$. Similarly, let $\psi\colon F\to F$ be the bijective $1$-cocycle with respect to $\vartheta$ induced by the rule
		$\psi(x)=x$ for $x\in X$. Denote by $\eta$ the congruence on $F$ generated by pairs $(xy,\lambda_x(y)\rho_y(x))$ for all $x,y\in X$.
		Clearly, we have $F/\eta\cong M$. Moreover, it follows from equation \eqref{YB1} that $\eta\s\Ker\vartheta$. Now, fix $x,y,z\in X$ and put
		\[u=\lambda_z(x)\in X\qquad\text{and}\qquad v=\lambda_{\rho_x(z)}(y)\in X.\]
		Then using equation \eqref{YB1} we get
		\begin{align*}
			\vartheta(z)(\psi(xy))
			& =\lambda_z(x\lambda_x(y))\\
			& =\lambda_z(x)\lambda_z(\lambda_x(y))\\
			& =\lambda_z(x)\lambda_{\lambda_z(x)}(\lambda_{\rho_x(z)}(y))\\
			& =u\lambda_u(v)\\
			& =\psi(uv).
		\end{align*}
		Furthermore, equations \eqref{YB1} and \eqref{YB2} yield
		\begin{align*}
			\vartheta(z)(\psi(\lambda_x(y)\rho_y(x)))
			& =\lambda_z(\lambda_x(y)\lambda_{\lambda_x(y)}(\rho_y(x)))\\
			&=\lambda_z(\lambda_x(y))\lambda_z(\lambda_{\lambda_x(y)}(\rho_y(x)))\\
			&=\lambda_{\lambda_z(x)}(\lambda_{\rho_x(z)}(y))\lambda_{\lambda_z(\lambda_x(y))}(\lambda_{\rho_{\lambda_x(y)}(z)}(\rho_y(x)))\\
			&=\lambda_{\lambda_z(x)}(\lambda_{\rho_x(z)}(y))\lambda_{\lambda_{\lambda_z(x)}(\lambda_{\rho_x(z)}(y))}
			(\rho_{\lambda_{\rho_x(z)}(y)}(\lambda_z(x)))\\
			&=\lambda_u(v)\lambda_{\lambda_u(v)}(\rho_v(u))\\
			&=\psi(\lambda_u(v)\rho_v(u)).
		\end{align*}
		Hence $\psi(\eta)=\{(\vartheta(z)(\psi(x)),\vartheta(z)(\psi(y))):(x,y)\in\eta\}$ for all $z\in F$. Concluding, the congruence
		$\eta$ satisfies both conditions (1) and (2) from Lemma \ref{lem:1}. Thus $\psi(\eta)$ is a congruence on $F$.
		Moreover, the congruence $\psi(\eta)$ is generated by pairs
		\[(\psi(xy),\psi(\lambda_x(y)\rho_y(x)))=( x\lambda_x(y),\lambda_x(y) \lambda_{\lambda_x(y)}(\rho_y(x)))\]
		for all $x,y\in X$. Therefore $F/\psi(\eta)\cong A$ and both statements (1) and (2) of the proposition are direct consequences of Lemmas
		\ref{lem:1} and \ref{lem:2}. Since statement (3) is obvious, the result is proved.
	\end{proof}
\end{prop}

It is worth to add that for a solution $(X,r)$ we can also define a ``right analog'' $A'(X,r)$ of the monoid $A(X,r)$ as
\[A'(X,r)=\free{X\mid \rho_y(x)y=\rho_{\rho_y(x)}(\lambda_x(y))\rho_y(x)\text{ for all }x,y\in X}.\]

If the solution $(X,r)$ is right non-degenerate then one can show (in a similar manner as in Proposition \ref{prop:1}) that
there exist a right action of $M(X,r)$ on $A'(X,r)$ and a bijective (right) $1$-cocycle $M(X,r)\to A'(X,r)$ with respect to this action.
Hence, one obtains that the structure monoid $M(X,r)$ is isomorphic to the regular submonoid $\{(\phi'(a),a):a\in A'(X,r)\}$ of the
semidirect product $\mc{G}'(X,r)^{\op}\ltimes A'(X,r)$, where $\mc{G}'(X,r)=\gr(\rho_x\mid x\in X)\s\Sym(X)$ and the map
$\phi'\colon A'(X,r)\to\mc{G}'(X,r)$ satisfies $\phi'(b)\phi'(a)=\phi'(\phi'(b)(a)b)$ for all $a,b\in A'(X,r)$.

\section{Structure of the monoid $A(X,r)$ and its algebra}\label{sec:2}

The following lemma and proposition are proved in \cite{So} and \cite{LYZ} for right non-degenerate solutions.
We include the proofs for completeness' sake.

\begin{lem}\label{lem:Sol1}
    Let $(X,r)$ be a left non-degenerate solution of the Yang--Baxter equation. Then
    \[\lambda_x\circ\sigma_y=\sigma_{\lambda_x(y)}\circ\lambda_x\] for all $x,y\in X$, where
    $\sigma_z$ for $z\in X$ is defined by the first equality in \eqref{sig}.
	\begin{proof}
		Let $z\in X$. By \eqref{sig} and \eqref{YB1} we get
		\[\lambda_x(\sigma_y(z))=\lambda_x(\lambda_y(\rho_{\lambda^{-1}_z(y)}(z)))=
		\lambda_{\lambda_x(y)}(\lambda_{\rho_y(x)}(\rho_{\lambda^{-1}_z(y)}(z))).\]
		Denote $t=\lambda_z^{-1}(y)\in X$. Then using \eqref{YB2} the previous becomes
		\[\lambda_{\lambda_{x}(y)}(\lambda_{\rho_{\lambda_{z}(t)}(x)}(\rho_t(z)))=
		\lambda_{\lambda_x(y)}(\rho_{\lambda_{\rho_z(x)}(t)}(\lambda_x(z)))=
		\lambda_{\lambda_x(y)}(\rho_{\lambda_{\rho_z(x)}(\lambda^{-1}_z(y))}(\lambda_x(z))).\]
		Applying \eqref{YB1} and \eqref{sig} once more, we get that this is equal to
		\[\lambda_{\lambda_x(y)}(\rho_{\lambda^{-1}_{\lambda_x(z)}(\lambda_x(y))}(\lambda_x(z)))=\sigma_{\lambda_x(y)}(\lambda_x(z)),\]
		and the result follows.
	\end{proof}
\end{lem}

\begin{prop}\label{prop:2}
    Assume that $(X,r)$ is a left non-degenerate solution of the Yang--Baxter equation. Define $s\colon X\times X\to X\times X$
    by $s(x,y)=(y,\sigma_y(x))$. Then $(X,s)$ is a left non-degenerate solution of the Yang--Baxter equation satisfying 
    $M(X,s)=A(X,s)=A(X,r)$. Moreover, the solution $(X,r)$ is bijective if and only if $(X,s)$ is bijective if and only if $(X,s)$ is right 
    non-degenerate, and in this case $(X,s)$ is called the rack solution associated to $(X,r)$.
    \begin{proof}
    Let $r_1={\id}\times r$, $r_2=r\times{\id}$, $s_1={\id}\times s$ and $s_2=s\times{\id}$. Our aim is to show that
    \[s_1\circ s_2\circ s_1=s_2\circ s_1\circ s_2.\] We shall prove that $s_i=J\circ r_i\circ J^{-1}$ for $1\le i\le 2$,
    where the map $J\colon X\times X\times X\to X\times X\times X$ is defined as $J(x,y,z)=(x,\lambda_x(y),\lambda_x(\lambda_y(z)))$,
    which clearly implies what we need. First note that $J$ is indeed a bijection with the inverse given by
    $J^{-1}(x,y,z)=(x,\lambda_x^{-1}(y),\lambda_{\lambda_x^{-1}(y)}^{-1}(\lambda_x^{-1}(z)))$.
    Put \[u=\lambda_x^{-1}(y)\in X\qquad\text{and}\qquad v=\lambda_x^{-1}(z)\in X.\]
    Then we have $J^{-1}(x,y,z)=(x,u,\lambda_u^{-1}(v))$. Since $\lambda_y(\rho_u(x))=\sigma_y(x)$ by \eqref{sig} and
    \[\lambda_y(\lambda_{\rho_u(x)}(\lambda_u^{-1}(v)))
        =\lambda_{\lambda_x(u)}(\lambda_{\rho_u(x)}(\lambda_u^{-1}(v)))
        =\lambda_x(\lambda_u(\lambda_u^{-1}(v)))
        =\lambda_x(v)=z\]
    by \eqref{YB1}, we obtain
    \begin{align*}
        (J\circ r_1\circ J^{-1})(x,y,z) & =J(r_1(x,u,\lambda_u^{-1}(v)))\\
        & =J(y,\rho_u(x),\lambda_u^{-1}(v))\\
        & =(y,\lambda_y(\rho_u(x)),\lambda_y(\lambda_{\rho_u(x)}(\lambda_u^{-1}(v))))\\
        & =(y,\sigma_y(x),z)\\
        & =s_1(x,y,z).
    \end{align*}
    Moreover, because $\lambda_v(\rho_{\lambda_u^{-1}(v)}(u))=\sigma_v(u)$ by \eqref{sig} and
    \[\lambda_x(\lambda_v(\rho_{\lambda_u^{-1}(v)}(u)))=\lambda_x(\sigma_v(u))=\sigma_{\lambda_x(v)}(\lambda_x(u))=\sigma_z(y)\]
    by Lemma \ref{lem:Sol1}, we conclude that
    \begin{align*}
        (J\circ r_2\circ J^{-1})(x,y,z) & =J(r_2(x,u,\lambda_{u}^{-1}(v)))\\
        & =J(x,v,\rho_{\lambda_u^{-1}(v)}(u))\\
        & =(x,z,\lambda_x(\lambda_v(\rho_{\lambda_u^{-1}(v)}(u))))\\
        & =(x,z,\sigma_z(y))\\
        & =s_2(x,y,z).
    \end{align*}
    Hence the first part of the result is proved. To finish the proof it is enough to observe that $s=I\circ r\circ I^{-1}$,
    where the bijection $I\colon X\times X\to X\times X$ is defined as $I(x,y)=(x,\lambda_x(y))$.
    \end{proof}
\end{prop}

Moreover, in case the solution $(X,r)$ is bijective, if we define $x\triangleleft y=\sigma_y(x)$ for $x,y\in X$ then the resulting structure 
$(X,\triangleleft)$ is a rack. If furthermore $\sigma_y(y)=y$ for all $y\in X$, this is a quandle (see also \cite{ESS}).

\begin{rem}\label{rem:1}
	Note that if $(X,r)$ is a bijective left non-degenerate solution then, by virtue of the defining relations,
	every element $x$ of $X$ is normal in $A=A(X,r)$. Hence each element of $A$ is normal, i.e., $aA=Aa$ for all $a\in A$.
	If $X=\{x_1,\dotsc,x_n\}$ is a finite set, then \[A=\{x_1^{k_1}\dotsm x_n^{k_n}:k_1,\dotsc,k_n\ge 0\}.\]
\end{rem}

In essence, Proposition \ref{prop:2} boils down to the following equality
\[\sigma_x\circ\sigma_y=\sigma_{\sigma_x(y)}\circ\sigma_x\] for all $x,y\in X$.
Moreover, the above equality assures that the action of $A=A(X,r)$ on $A$, given as
\[\sigma_a=\sigma_{a_r}\circ\dotsb\circ\sigma_{a_1}\qquad\text{and}\qquad\sigma_a(b)=\sigma_a(b_1)\dotsm\sigma_a(b_s)\]
for $a=a_1\dotsm a_r\in A$ and $b=b_1\dotsm b_s\in A$, where $a_i,b_j\in X$, is well-defined.
We shall freely use this fact throughout the paper.

\begin{prop}
	Assume that $(X,r)$ is a bijective left non-degenerate solution of the Yang--Baxter equation.
	Then there exist a set $I$ and $\sigma$-invariant submonoids $A_i$ of $A=A(X,r)$ for $i\in I$ 
	(i.e., $\sigma_a(A_i)\s A_i$ for all $a\in A$ and $i\in I$) such that $A$ is the subdirect product of the
	family $(A_i)_{i\in I}$ and $A_i A_j = A_j A_i$ for all $i,j\in I$. Furthermore, if $X$ is a finite set, then $I$ can be taken as a finite set.
	\begin{proof}
	    For $x,y\in X$ we declare that $x\sim y$ if and only if there exists $a\in A$ such that
		$\sigma_a(x)=y$. It is clear that $\sim$ is an equivalence relation on $X$. So, let $X=\bigcup_{i\in I}X_i$ be the partition of $X$
		with respect to $\sim$. Let $A_i=\free{X_i}$ for $i\in I$ denote the submonoid of $A$ generated by $X_i$. Clearly, each monoid
		$A_i$ is $\sigma$-invariant. Moreover, as each element of $A$ is normal, it follows that $A$ is the subdirect product of
		the family $(A_i)_{i\in I}$.
	\end{proof}
\end{prop}

\begin{lem}\label{lem:3}
	Assume that $(X,r)$ is a finite bijective left non-degenerate solution of the Yang--Baxter equation.
	Then there exists $d\ge 1$ such that $a^d$ is a central element of $A=A(X,r)$ for each $a\in A$.
	\begin{proof}
		As $X$ is a finite set, it follows that there exists $d\ge 1$ (we can choose $d$ as a divisor of $n!$,
		where $n=|X|$) such that $\sigma_a^d=\id$ for each $a\in A$. Now, if $b\in A$ then
		\[ba^d=a\sigma_a(b)a^{d-1}=\dotsb=a^d\sigma_a^d(b)=a^db,\] and the result follows.
	\end{proof}
\end{lem}

\begin{rem}
    Moreover, if $ac=bc$ or $ca=cb$ holds for some $a,b,c\in A$ then $az^i=bz^i$ for some $i\ge 1$, where
	$z=x_1^d\dotsm x_n^d\in\Z(A)$ (here $d\ge 1$ is defined as in Lemma \ref{lem:3}). Hence the monoid $A$
	is left cancellative if and only if it is right cancellative if and only if the central elements of $A$ are cancellable.
\end{rem}

\begin{thm}\label{thm:1}
	Assume that $(X,r)$ is a finite bijective left non-degenerate solution of the Yang--Baxter equation. Then $A=A(X,r)$ is a central-by-finite
	monoid, i.e., $A=\bigcup_{f\in F}Cf$ for a central submonoid $C\s A$ and a finite subset $F\s A$. In particular, if $K$ is a field then $K[A]$
	is a finite module over a central affine subalgebra of $K[A]$. Hence, $K[A]$ is a Noetherian PI-algebra satisfying
	\[\clK K[A]=\GK K[A]=\rk A\le|X|\] and the equality holds if and only if the solution $(X,r)$ is involutive.
	\begin{proof}
		Write $X=\{x_1,\dotsc,x_n\}$ with $n=|X|$. By Lemma \ref{lem:3} we know that there exists $d\ge 1$
		such that $x_1^d,\dotsc,x_n^d$ are central elements of $A$. Define $C=\free{x_1^d,\dotsc,x_n^d}$.
		Clearly $C$ is a central submonoid of $A$. Moreover, Remark \ref{rem:1} yields $A=\bigcup_{f\in F}Cf$,
		where \[F=\{x_1^{k_1}\dotsm x_n^{k_n}:0\le k_1,\dotsc,k_n<d\}\s A.\] In particular,
		$K[A]=\sum_{f\in F}K[C]f$ is a finite module over the central affine subalgebra $K[C]$ of $K[A]$.
		Therefore, the algebra $K[A]$ is Noetherian and PI. Hence, a result of Anan'in \cite{An} implies that $A$
		is a linear monoid and then \cite[Proposition 1, p.~221, Proposition 7, p.~280--281, and Theorem 14, p.~284]{O1} yield
		$\clK K[A]=\GK K[A]=\rk A$. Because the commutative algebra $K[C]$ can be generated by $n$ elements,
		we get $\GK K[A]=\GK K[C]\le n$, as desired.
		
		Finally, it is clear that if $(X,r)$ is involutive then the equality $\clK K[A]=\GK K[A]=\rk A=n$ holds
		as $A$ is a free abelian monoid of rank $n$. Whereas if $(X,r)$ is not involutive then we claim that
		$\sigma_x(y)\ne y$ for some $x,y\in X$. Indeed, otherwise $\sigma_x=\id$ for all $x\in X$ and then
		$\lambda_x=\hat{\lambda}_x$ by \eqref{sig}. Since
		\[\rho_y(x)=\hat{\lambda}^{-1}_{\lambda_x(y)}(x)=\lambda^{-1}_{\hat{\lambda}_x(y)}(x)=\hat{\rho}_y(x),\]
		we get $\rho_y=\hat{\rho}_y$ for all $y\in X$. Thus $r=r^{-1}$ and $(X,r)$ is involutive, a contradiction.
		Hence $\sigma_x(y)\ne y$ for some $x,y\in X$. Now $y^dx=x\sigma_x(y)^d=\sigma_x(y)^dx$ yields $(y^d-\sigma_x(y)^d)x^d=0\in K[C]$.
		Therefore, if $\mf{p}$ is a prime ideal of $K[C]$ then $x^d\in\mf{p}$ or $y^d-\sigma_x(y)^d\in\mf{p}$.
		Thus the commutative algebra $K[C]/\mf{p}$ can be generated by less than $n$ elements and it follows that $\clK K[A]=\clK K[C]<n$.
	\end{proof}
\end{thm}

For definitions of all homological notions used below we refer to \cite{BG,Ga2,Lav,SZ}.

\begin{thm}\label{prop:4}
	Assume that $(X,r)$ is a finite bijective left non-degenerate solution of the Yang--Baxter equation.
	Let $A=A(X,r)$. If $K$ is a field then the following conditions are equivalent:
	\begin{enumerate}
		\item $(X,r)$ is an involutive solution.
		\item $A$ is a free abelian monoid of rank $|X|$.
		\item $A$ is a cancellative monoid.
	    \item $\rk A=|X|$.
		\item $K[A]$ is a prime algebra.
		\item $K[A]$ is a domain.
	    \item $\clK K[A]=|X|$.
	    \item $\GK K[A]=|X|$.
	    \item $\id K[A]=|X|$.
	    \item $\gld K[A]=|X|$.
	    \item $K[A]$ is an Auslander--Gorenstein algebra.
	    \item $K[A]$ is an Auslander-regular algebra.
	\end{enumerate}
	\begin{proof}
		Clearly $(1)\Longrightarrow(2)\Longrightarrow(10)$ and $(2)\Longrightarrow(6)\Longrightarrow(5)$ and $(12)\Longrightarrow(11)$.
		Moreover, we have $(1)\Longleftrightarrow(4)\Longleftrightarrow(7)\Longleftrightarrow(8)$ by Theorem \ref{thm:1} and
		$(5)\Longrightarrow(3)$ by Remark \ref{rem:1}. Further, $(10)\Longrightarrow(12)$ and $(11)\Longleftrightarrow(9)\Longrightarrow(8)$
		follows by Theorem \ref{thm:1} and \cite[Theorem 1, p.~126]{BG} (see also \cite{SZ}). Summarizing, we have the following diagram of implications
		
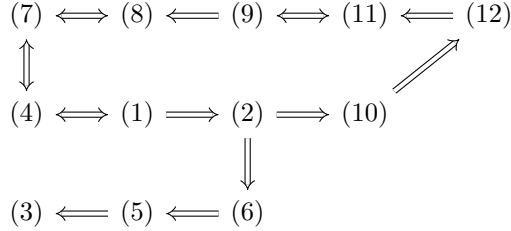
\[\begin{tikzcd}[row sep=20pt,column sep=20pt]
		    (7)\ar[d,Leftrightarrow] & (8)\ar[l,Leftrightarrow] & (9)\ar[l,Rightarrow] & (11)\ar[l,Leftrightarrow] & (12)\ar[l,Rightarrow]\\
		    (4)\ar[r,Leftrightarrow] & (1)\ar[r,Rightarrow] & (2)\ar[r,Rightarrow]\ar[d,Rightarrow] & (10)\ar[ru,Rightarrow]\\
		    (3) & (5)\ar[l,Rightarrow] & (6)\ar[l,Rightarrow]
		\end{tikzcd}\]
		Therefore, it is enough to show that $(3)\Longrightarrow(2)$.
		So assume $(3)$ and observe first that then $xx=x\sigma_x(x)$ yields $\sigma_x(x)=x$
		for each $x\in X$. Now, choose $d\ge 1$ such that $a^d\in\Z(A)$ for each $a\in A$. Then for $x,y\in X$ we have 
		\[y^dx=y^{d-1}x\sigma_x(y)=\dotsb=x\sigma_x(y)^d=\sigma_x(y)^dx,\] which leads to $\sigma_x(y)^d=y^d$.
		Because the elements $y^d$ and $\sigma_x(y)^d$ cannot be rewritten using the defining relations of $A$
		(the only way to rewrite the word $z^d$ for $z\in X$ would be to use a relation of the form
		$z\sigma_z(z)=zz=\sigma_z^{-1}(z)z$), we conclude that $\sigma_x(y)=y$. Hence $\sigma_x=\id$ and (2) follows. This finishes the proof.
	\end{proof}
\end{thm}

Assume that $(X,r)$ is a finite bijective left non-degenerate solution of the Yang--Baxter equation. Let $A=A(X,r)$ and define
\[\eta_A=\{(a,b)\in A\times A:ac=bc\text{ for some }c\in A\}.\] By Remark \ref{rem:1} it follows that $\eta_A$
is the cancellative congruence of the monoid $A$, that is the smallest congruence $\eta$ on $A$ such that the quotient monoid
$A/\eta$ is cancellative. Moreover, $\eta_A=\bigcup_{i=1}^\infty\eta_i$ is a union of the ascending chain of congruences 
$\eta_i=\{(a,b)\in A\times A:az^i=bz^i\}$ (here $z=\prod_{x\in X}x^d\in\Z(A)$ is defined as in Remark \ref{rem:1}).
Note that the lattice of congruences on $A$ can be embedded into the lattice of ideals of the algebra $K[A]$ (here $K$ is an arbitrary field), 
by associating to a congruence $\eta$ on $A$ the ideal \[I(\eta)=\Span_K\{a-b:(a,b)\in \eta\},\] the $K$-linear span of the set consisting
of all elements $a-b$ with $(a,b)\in\eta$. We conclude by Theorem \ref{thm:1} that the monoid $A$ satisfies the ascending chain condition
on congruences. Hence there exists $t\ge 1$ such that $\eta_i=\eta_t$ for each $i\ge t$, and thus $\eta_A=\eta_t$.
Therefore, we have proved the following result.

\begin{prop}\label{prop:5}
	Assume that $(X,r)$ is a finite bijective left non-degenerate solution of the Yang--Baxter equation.
	Then there exists $t\ge 1$ such that \[\eta_A=\{(a,b)\in A\times A:az^i=bz^i\}\] for all $i\ge t$,
	where $z=\prod_{x\in X}x^d\in\Z(A)$ is defined as in Remark \ref{rem:1}.
	In particular, the ideal $Az^t$ is cancellative and if $K$ is a field then $I(\eta_A)=\Ann_{K[A]}(z^i)$ for all $i\ge t$.
\end{prop}

\section{Prime ideals of $A(X,r)$ and $K[A(X,r)]$}\label{sec:3}

We shall begin this section with the following description of prime ideals of the monoid $A=A(X,r)$.
The prime spectra of $A$ and $K[A]$ (for a field $K$) are denoted as $\Spec(A)$ and $\Spec(K[A])$, respectively.

\begin{prop}\label{prop:6}
	Assume that $(X,r)$ is a finite bijective left non-degenerate solution of the Yang--Baxter equation.
	Let $A=A(X,r)$ and \[\mc{Z}=\mc{Z}(X,r)=\{Z\s X:\vn\ne Z\ne X\text{ and }\sigma_x(Z)=Z\text{ for all }x\in X\setminus Z\}.\]
	Define $P(Z)=\bigcup_{z\in Z}Az$ for $Z\in\mc{Z}$. Then the maps
	\[\mc{Z}\to\Spec(A)\colon Z\mapsto P(Z)\qquad\text{and}\qquad\Spec(A)\to\mc{Z}\colon P\mapsto X\cap P\]
	are mutually inverse bijections.
	\begin{proof}
		Since the elements of $A$ are normal it is clear that if $P\in\Spec(A)$ then $\vn\ne X\cap P\ne X$ and
		$P=\bigcup_{x\in X\cap P}Ax$. Moreover, if $x\in X\cap P$ and $y\in X\setminus P$ then $y\sigma_y(x)=xy\in P$.
		Hence $y\notin P$ implies $\sigma_y(x)\in P$. Therefore, $\sigma_y(X\cap P)=X\cap P$ and $X\cap P\in\mc{Z}$.
		
		Conversely, if $Z\in\mc{Z}$ then we claim that $P(Z)$ is a prime ideal of $A$. To show this observe that if
		$x_1\dotsm x_n\in P(Z)$ for some $x_1,\dotsc,x_n\in X$ then necessarily $x_i\in P(Z)$ for some $1\le i\le n$.
		Otherwise $x_1,\dotsc,x_n\in X\setminus Z$ and then each word in the free monoid on $X$ representing the element
		$x_1\dotsm x_n\in A$ must be a product of letters in $X\setminus Z$, which leads to a contradiction. Indeed, if
		$x,y\in X\setminus Z$ then the only way to rewrite the word $xy$ is to use one of the relations $xy=y\sigma_y(x)$
		and $xy=\sigma_x^{-1}(y)x$. Since $\sigma_x(Z)=\sigma_y(Z)=Z$, we get
		$\sigma_x(X\setminus Z)=\sigma_y(X\setminus Z)=X\setminus Z$.
		Hence both $\sigma_y(x)$ and $\sigma_x^{-1}(y)$ are elements of $X\setminus Z$.
	\end{proof}
\end{prop}

Our next result provides an inductive description of all prime ideals of the monoid algebra $K[A(X,r)]$ over a field $K$ in terms
of prime ideals of group algebras over $K$ of certain finitely generated FC-groups (finite conjugacy groups) closely related
to the monoid $A(X,r)$. Recall that for such a group $G=\Delta(G)$ the torsion elements form a finite characteristic subgroup
$G^+=\Delta^+(G)$ such that $G/G^+$ is a finitely generated free abelian group (see, e.g., \cite[Section 4.1]{P1}).

\begin{prop}\label{prop:7}
	Assume that $(X,r)$ is a finite bijective left non-degenerate solution of the Yang--Baxter equation. Let $A=A(X,r)$ and
	$\mc{Z}=\mc{Z}(X,r)$. If $K$ is a field and $P$ is a prime ideal of the algebra $K[A]$ then $X\cap P\in\mc{Z}\cup\{\vn,X\}$.
	Moreover, for such a prime ideal $P$ the following properties hold:
	\begin{enumerate}
		\item there exists an inclusion preserving bijection between the set of prime ideals $Q$ of $K[A]$ with
		the property $X\cap Q=X\cap P$ and the set of all prime ideals of the algebra $K[A\setminus P]$.
		Moreover, the monoid $A\setminus P$ has the following presentation
		\[A\setminus P\cong\free{X\setminus P\mid xy=y\sigma_y(x)\text{ for all }x,y\in X\setminus P}.\]
		\item there exists an inclusion preserving bijection between the set of prime ideals $Q$ of $K[A]$ satisfying
		$Q\cap A=\vn$ and the set of all prime ideals of the group algebra $K[G]$, where
		\[G=\gr(X\mid xy=y\sigma_y(x)\text{ for all }x,y\in X).\]
		Furthermore, the cancellative monoid $\ov{A}=A/\eta_A$ has a group of quotients, which is  equal to
		the central localization $\ov{A}\free{z}^{-1}$ for some $z\in\Z(\ov{A})$, and $G\cong\ov{A}\free{z}^{-1}$.
		Clearly, $G$ is a finitely generated FC-group.
	\end{enumerate}
	\begin{proof}
		Clearly $P\cap A$ is a prime ideal of $A$ and $X \cap (P\cap A) = X \cap P$.
		Hence, from Proposition \ref{prop:6}, we get that $X \cap P\in\mc{Z}$ if $\vn \ne X \cap P \ne X$.
		Therefore, the first part of the proposition follows.
		
		Since $Q\cap A=\bigcup_{x\in Q\cap X}Ax$, it is clear that the set of prime ideals $Q$ of $K[A]$ with the
		property $X\cap Q=X\cap P$ is in an inclusion preserving bijection with the set of all prime ideals of the algebra
		$K[A]/K[P\cap A]\cong K_0[A/(P\cap A)]$, the contracted semigroup algebra of $A/(P\cap A)$ (recall that the contracted
		semigroup algebra $K_0[S]$, for a semigroup $S$ with zero element $\theta$, is defined as $K[S]/K\theta$).
		By Proposition \ref{prop:6} we get $A\setminus P=\free{X\setminus P}\s A$ and also
		$A\setminus P\cong\free{X\setminus P\mid xy=y\sigma_y(x)\text{ for all }x,y\in X\setminus P}$.
		Therefore, $A/(P\cap A)\cong(A\setminus P)\cup\{\theta\}$.
		Hence $K_0[A/(P\cap A)]\cong K[A\setminus P]$, and (1) follows.
		
		Assume now that $Q$ is a prime ideal of $K[A]$ such that $Q\cap A=\vn$. We claim that $Q$ contains the ideal
		$I(\eta_A)$. Indeed, if $a,b\in A$ satisfy $ac=bc$ for some central element $c\in A$ then
		\[(a-b)K[A]c=(a-b)cK[A]=0\s Q.\] Since $c\notin Q$, we get $a-b\in Q$. Therefore the ideals of $K[A]$ intersecting $A$
		trivially correspond bijectively to the prime ideals of the algebra $K[A]/I(\eta_A)\cong K[\ov{A}]$,
		and hence also to the prime ideals of the central localization $K[\ov{A}]\free{z}^{-1}\cong K[\ov{A}\free{z}^{-1}]$,
		where $z=\prod_{x\in X}x^d\in\Z(A)$ (here $d\ge 1$ is defined as in Lemma \ref{lem:3}).
		So, it remains to show that the group $G$ is isomorphic to $\ov{A}\free{z}^{-1}$, which is clearly equal
		to the group of quotients of the monoid $\ov{A}$. Observe that the natural morphism $A\to G$ factors
		through $\ov{A}$ and thus also through $\ov{A}\free{z}^{-1}$. Hence we get a natural morphism of groups
		$\varphi\colon\ov{A}\free{z}^{-1}\to G$.
		
		Surjectivity of $\varphi$ is obvious. Let $\ov{a},\ov{b}\in \ov{A}\free{z}^{-1}$ such that $\varphi(\ov{a})=\varphi(\ov{b})$.
		We may assume that $\ov{a},\ov{b}\in\ov{A}$ as we can multiply them by their highest denominator in $z$. Consider $a$ and $b$ as words
		in the free group $F$ with generators in $X$. By adding the relation $xy = y \sigma_y(x)$ on the free group $F$,
		we get the group $G$, where the words corresponding to $a$ and $b$ are equal. Hence, they are equal in every
		group generated in $X$, satisfying the relation $xy=y\sigma_y(x)$. Thus, $\ov{a}=\ov{b}$ in the group
		$\ov{A}\free{z}^{-1}$. Therefore, $\varphi$ is injective, which finishes the proof.
	\end{proof}
\end{prop}

As Example \ref{ex:2} shows, it is possible that the algebra $K[A(X,r)]$ does not admit prime ideals
intersecting the monoid $A(X,r)$ non-trivially, even if the group $\Sigma=\gr(\sigma_x\mid x\in X)\s\Sym(X)$ is cyclic.

However, it is clear that for a prime ideal $P$ of $K[A(X,r)]$ we have $P\cap A(X,r)\ne\vn$ if and only if $z\in P$,
where the element $z=\prod_{x\in X}x^d\in\Z(A(X,r))$ is defined as in Remark \ref{rem:1}.
Thus, the maximal ideal containing $1-z$ is a prime ideal that intersects $A(X,r)$ trivially.
Hence the algebra $K[A(X,r)]$ always has minimal prime ideals intersecting the monoid $A(X,r)$ trivially.

\begin{ex}\label{ex:2}
	Let $X$ be a finite non-empty set. Fix $\sigma\in\Sym(X)$ and define $r\colon X\times X\to X\times X$ as
	$r(x,y)=(y,\sigma(x))$. Clearly $(X,r)$ is a bijective non-degenerate solution of the Yang--Baxter equation.
	Let $A=A(X,r)$ and $K$ be a field. If $x,y\in X$ then $y(x-\sigma(x))=(x-\sigma(x))\sigma(y)$
	in $K[A]$. Thus $x-\sigma(x)$ is a normal element of $K[A]$ and hence $K[A](x-\sigma(x))$ is an ideal of $K[A]$.
	We claim that this ideal is nilpotent. Indeed, observe first that the equality $x\sigma(x)=\sigma(x)^2$ yields
	\[(x-\sigma(x))^2=x^2-x\sigma(x)-\sigma(x)x+\sigma(x)^2=(x-\sigma(x))x,\] which leads to
	$(x-\sigma(x))^{n+1}=(x-\sigma(x))x^n$ for each $n\ge 1$. In particular, if $d\ge 1$ is equal to the order
	$\sigma$, then $(x-\sigma(x))^d=0$. Indeed, if $d=1$ then the equality is obvious. Whereas, if $d\ge 2$ then
	\[(x-\sigma(x))^d=(x-\sigma(x))x^{d-1}=x^d-\sigma(x)x^{d-1}=x^d-x^{d-1}\sigma^d(x)=0.\]
	Thus $(K[A](x-\sigma(x)))^d=K[A](x-\sigma(x))^d=0$, as claimed. Hence the ideal $P=\sum_{x\in X}K[A](x-\sigma(x))$
	is nilpotent. Note that if $x,y\in X$ then $xy-yx=x(y-\sigma(y))\in P$. Moreover, if $x\in X$ and $n\ge 1$ then 
	\[x-\sigma^n(x)=\sum_{i=1}^n(\sigma^{i-1}(x)-\sigma^i(x))\in P.\] These facts easily lead to a conclusion that
	$K[A]/P\cong K[t_1,\dotsc,t_s]$, the commutative polynomial algebra in $s$ commuting variables, where $s$ is the number
	of disjoint cycles in the decomposition of $\sigma$. Hence the ideal $P$ is also semiprimitive. Therefore,
	$\mc{B}(K[A])=\mc{J}(K[A])=P$ is the unique minimal prime ideal of $K[A]$. In particular, $\clK K[A]=s$
	may be equal to any prescribed integer between $1$ and $|X|$.
\end{ex}

Moreover, as Example \ref{ex:3} shows, the description of the minimal primes of the algebra $K[A(X,r)]$
depends on the characteristic of a base field $K$.

\begin{ex}\label{ex:3}
	Consider the solution $(X,r)$ defined in Example \ref{ex:1}. Let $A=A(X,r)$ and assume that $K$ is a field.
	The following facts can be verified (using theory of Gr\"obner bases and the fact that the algebras under
	consideration are $\mb{Z}$-graded).
	If $\ch K=3$ then the minimal prime ideals of the algebra $K[A]$ are of the form:
	\[P_1=(x_2,x_3),\qquad P_2=(x_1,x_3),\qquad P_3=(x_1,x_2),\qquad P_4=(x_1-x_2,x_2-x_3).\]
	Whereas, if $\ch K\ne 3$ then the minimal primes of $K[A]$ consist of the ideals $P_1,P_2,P_3,P_4$ together with the ideal 
	\[P_5=(x_1+x_2+x_3,x_1^2-x_2^2,x_2^2-x_3^2)=(x_1+x_2+x_3,x_1^2-x_2^2).\]
	Furthermore, if $\ch K=3$ then $0\ne x_1(x_2-x_3)\in\mc{B}(K[A])$, whereas if $\ch K\ne 3$ then the algebra $K[A]$ is semiprime.
\end{ex}

Our next aim is to determine the classical Krull dimension (which is equal to the Gelfand--Kirillov dimension;
see Theorem \ref{thm:1}) of the algebra $K[A(X,r)]$ over a field $K$ in terms of certain purely combinatorial
properties of the permutations $\sigma_x$ for $x\in X$.

\begin{thm}\label{thm:2}
	Assume that $(X,r)$ is a finite bijective left non-degenerate solution of the Yang--Baxter equation.
	Let $A=A(X,r)$ and $G=\gr(X \mid xy = y\sigma_y(x)\text{ for all }x,y\in X)$.
	If $K$ is a field and $P$ is a minimal prime ideal of the algebra $K[A]$ satisfying $P\cap A=\vn$
	then \[\clK K[A]/P=\clK K[G]=s,\] where $s$ is the number of orbits of $X$ with respect to the action of the group
	$\Sigma=\gr(\sigma_x\mid x\in X)\s\Sym(X)$.
	\begin{proof}
		By Lemma \ref{lem:3} there exists $d\ge 1$ such that $\sigma_x^d=\id$ for each $x\in X$.
		Define $C=\free{x^d\mid x\in X}\s A$. Clearly, $C$ is a central submonoid of $A$. First, we shall prove that
		\begin{equation}
			\clK K[A]/P\le s.\label{clK1}
		\end{equation}
		Note that if $x\in X$ and $a\in A$ then $ax^d=x^d\sigma_x^d(a)=x^da=a\sigma_a(x)^d$. Hence
		\[aK[A](x^d-\sigma_a(x)^d)=K[A]a(x^d-\sigma_a(x)^d)=0\s P.\] Thus $a\notin P$ leads to $x^d-\sigma_a(x)^d\in P$.
		Since $K[C]$ is a central subalgebra of $K[A]$, the ideal $\mf{p}=P\cap K[C]$ is prime. Moreover,
		$x^d-\sigma_a(x)^d\in\mf{p}$ for all $x\in X$ and $a\in A$. Therefore, $\clK K[C]/\mf{p}\le s$ because the commutative
		algebra $K[C]/\mf{p}$ can be generated by $s$ elements (the image of the set $\{x^d:x\in X\}\s C$
		in $K[C]/\mf{p}$ has cardinality $\le s$). Since $K[A]/P$ is PI-algebra, which is a finite module over the
		central subalgebra $K[C]/\mf{p}$, we conclude by \cite[Theorem 13.8.14]{MR} that $\clK K[A]/P=\clK C/\mf{p}\le s$, as desired.
		
		Next we shall prove that
		\begin{equation}
			\clK K[G]\le\clK K[A]/P.\label{clK2}
		\end{equation}
		By Proposition \ref{prop:7} we know that $P$ corresponds bijectively to a minimal prime ideal $\ov{P}$ of the algebra
		$K[A/\eta_A]$ and also to a minimal prime ideal $P_G$ of the group algebra $K[G]$. Let $z=\prod_{x\in X}x^d\in C$.
		The images of $z$ in the algebras $K[A/\eta_A]$ and $K[A]/P$, still denoted by $z$, are central and regular elements
		of these algebras. Since $K[A/\eta_A]\free{z}^{-1}\cong K[G]$ by Proposition \ref{prop:7}, it follows easily that
		\[(K[A]/P)\free{z}^{-1}\cong(K[A/\eta_A]/\ov{P})\free{z}^{-1}\cong K[G]/P_G.\]
		Hence \[\clK K[G]/P_G=\clK(K[A]/P)\free{z}^{-1}\le\clK K[A]/P.\] Thus, to prove \eqref{clK2}
		it is enough to show that $\clK K[G]/P_G=\clK K[G]$.
		By \cite[Lemma 4.1.8]{P1} we know that there exists a finitely generated free abelian subgroup $F\s G$ of finite index.
		Therefore, $K[G]\cong K[F]*(G/F)$, a crossed product of the finite group $G/F$ over the Laurent polynomial algebra $K[F]$.
		Hence \cite[Theorem 16.6]{P2} guarantees that \[\h Q=\h Q\cap K[F]\qquad\text{and}\qquad\clK K[G]/Q=\clK K[F]/(Q\cap K[F])\]
		for each $Q\in\Spec(K[G])$. By Schelter's theorem (see \cite[Theorem 13.10.12]{MR}) we get
		\[\h\mf{p}+\clK K[F]/\mf{p}=\clK K[F]\] for each $\mf{p}\in\Spec(K[F])$.
		Since $\clK K[G]=\clK K[F]$, we conclude that \[\h Q+\clK K[G]/Q=\clK K[G]\] for each $Q\in\Spec(K[G])$.
		In particular, as $\h P_G=0$, we obtain $\clK K[G]/P_G=\clK K[G]$, as desired.
		
		Finally, let us observe that the ideal $P_0$ of $K[G]$ generated by elements $x-y$ for all $x,y\in X$ which are in the same orbit
		of $X$ with respect to the action of $\Sigma$ satisfies $K[G]/P_0\cong K[t_1^{\pm 1},\dotsc,t_s^{\pm 1}]$, the Laurent polynomial
		algebra in $s$ commuting variables. Hence
		\begin{equation}
			s=\clK K[G]/P_0\le\clK K[G].\label{clK3}
		\end{equation}
		Putting \eqref{clK1}, \eqref{clK2} and \eqref{clK3} together we get $\clK K[A]/P=\clK K[G]=s$, which finishes the proof.
	\end{proof}
\end{thm}

Motivated by Propositions \ref{prop:6}, \ref{prop:7} and Theorem \ref{thm:2} we define \[\Sigma_Z=\gr(\sigma_x\mid x\in X\setminus Z)\s\Sym(X)\]
and \[s(Z)=\text{the number of orbits of }X\setminus Z\text{ with respect to the action of }\Sigma_Z\] for each $Z\in\mc{Z}_0=\mc{Z}\cup\{\vn\}$,
where $\mc{Z}=\mc{Z}(X,r)$.

By Proposition \ref{prop:6} we know that all sets in $\mc{Z}$ are of the form $X\cap P$ for a prime ideal $P$
of $A=A(X,r)$. On the other hand, Proposition \ref{prop:7} assures that if $Q$ is a prime ideal of the algebra $K[A]$
over a field $K$ then $X\cap Q\in\mc{Z}_0$ or $X\cap Q=X$. But if $Q$ is a minimal prime ideal of $K[A]$ then the latter
possibility is excluded. Indeed, otherwise $Q$ would strictly contain the prime ideal $Q_0$ generated by elements
$x-y$ for all $x,y\in X$. However, as Example \ref{ex:4} shows, not all sets in $\mc{Z}_0$ are of the form
$X\cap Q$ for a minimal prime ideal $Q$ of $K[A]$.

\begin{ex}\label{ex:4}
    Let $X=\{x_1,x_2,x_3,x_4,x_5\}$. Define $\sigma_1=\sigma_2=(1,2)$, $\sigma_3=\sigma_5=\id$ and $\sigma_4=(3,5)$.
    Let $r\colon X\times X\to X\times X$ be defined as $r(x_i,x_j)=(x_j,x_{\sigma_j(i)})$. It is easy to check that
    $(X,r)$ is a bijective non-degenerate solution of the Yang--Baxter equation. If $K$ is a field and $A=A(X,r)$ then
    \[x_4(x_3-x_5)=0,\qquad x_1(x_1-x_2)=0,\qquad x_2(x_1-x_2)=0\] in $K[A]$.
    Clearly, the first equality assures that each prime ideal $P$ of $K[A]$ contains $x_4$ or $x_3-x_5$.
    Moreover, the second and third equalities guarantee that $x_1-x_2\in P$. Because $P_1=(x_1-x_2,x_4)$ and $P_2=(x_1-x_2,x_3-x_5)$
    are prime ideals of $K[A]$ (actually, we have $K[A]/P_1\cong K[A]/P_2\cong K[t_1,t_2,t_3]$, the polynomial algebra in
    three commuting variables), $P_1$ and $P_2$ are the only minimal prime ideals of $K[A]$.
    However, the set $Z=\{x_3,x_4\}\in\mc{Z}(X,r)$ satisfies $Z\ne X\cap P_1=\{x_4\}$ and $Z\ne X\cap P_2=\vn$.
\end{ex}

Note that Example \ref{ex:4} shows also that the algebra $K[A]$, where $A=A(X,r)$, may contain minimal prime ideals
of mixed type (i.e., prime ideals $P$ of $K[A]$ satisfying $P\cap A\ne\vn$ but $P\ne K[P\cap A]$),
even if the group $\Sigma=\gr(\sigma_x \mid x\in X)\s\Sym(X)$ is abelian. This is in contrast to what
happens in the cancellative case (see \cite{JO}).

Moreover, Example \ref{ex:5} shows that it is possible that the algebra $K[A]$ contains prime
ideals of the form $P=K[P\cap A]$, even if each orbit of $X$ with respect to the action of the group
$\Sigma$ has cardinality larger than $1$.

\begin{ex}\label{ex:5}
	Let $X=\{x_1,x_2,x_3,x_4\}$. Define $\sigma_1=\sigma_2=\id$ and $\sigma_3=\sigma_4=(1,2)(3,4)$.
    Moreover, let $r\colon X\times X\to X\times X$ be defined as $r(x_i,x_j)=(x_j,x_{\sigma_j(i)})$,
    a bijective non-degenerate solution of the Yang--Baxter equation.
    If $K$ is a field and $A=A(X,r)$ then
    \[x_3(x_3-x_4)=0,\qquad x_4(x_3-x_4)=0,\qquad x_3(x_1-x_2)=0,\qquad x_4(x_1-x_2)=0\] in $K[A]$.
    The above equalities assure that each prime ideal of $K[A]$ contains $x_3,x_4$ or $x_1-x_2,x_3-x_4$.
    Since $P=(x_3,x_4)$ is a prime ideal of $K[A]$ (actually, we have $K[A]/P\cong K[t_1,t_2]$,
    the polynomial algebra in two commuting variables), it is a minimal prime ideal of $K[A]$.
\end{ex}

\begin{thm}
	Assume that $(X,r)$ is a finite bijective left non-degenerate solution of the Yang--Baxter equation.
	Let $A=A(X,r)$ and $\mc{Z}_0=\mc{Z}\cup\{\vn\}$, where $\mc{Z}=\mc{Z}(X,r)$. If $K$ is a field then
	\[\clK K[A]=\max\{s(Z):Z\in\mc{Z}_0\}.\]
	\begin{proof}
		Define $s=\max\{s(Z):Z\in\mc{Z}_0\}$. If $P$ is a minimal prime ideal of $K[A]$ then $X\cap P\in\mc{Z}_0$
		and $P$ corresponds to a minimal prime ideal $P'$ of the algebra $K[A\setminus P]$ such that
		$P'\cap(A\setminus P)=\vn$ and $K[A]/P\cong K[A\setminus P]/P'$ (see Proposition \ref{prop:7}).
		Therefore, Theorem \ref{thm:2} implies that \[\clK K[A]/P=\clK K[A\setminus P]/P'=s(X\cap P)\le s.\]
		Since we have $\clK K[A]=\clK K[A]/P$ for some minimal prime ideal $P$ of $K[A]$, the inequality $\clK K[A]\le s$ follows.
		
		To show that $\clK K[A]\ge s$ we have to check that $\clK K[A]\ge s(Z)$ for each $Z\in\mc{Z}_0$.
		So, let us fix $Z\in\mc{Z}_0$. If $Z=\vn$ then we are done by Theorem \ref{thm:2}.
		Whereas, if $Z\in\mc{Z}$ then $A=P(Z)\cup A(Z)$, where $P(Z)=\bigcup_{z\in Z}Az$ and
		$A(Z)=A\setminus P(Z)=\free{X\setminus Z}\s A$ is the submonoid of $A$ generated by $X\setminus Z$.
		Therefore, $K[A]/K[P(Z)]\cong K_0[A/P(Z)]\cong K[A(Z)]$, which leads to
		\[\clK K[A]\ge\clK K[A]/K[P(Z)]=\clK K[A(Z)]\ge s(Z),\]
		where the last inequality is a consequence of the fact that the ideal $P_0$ of $K[A(Z)]$, generated
		by elements of the form $x-y$ for all $x,y\in X\setminus Z$ which are in the same orbit of $X\setminus Z$
		with respect to the action of the group $\Sigma_Z$, satisfies $K[A(Z)]/P_0\cong K[t_1,\dotsc,t_{s(Z)}]$,
		the polynomial algebra in $s(Z)$ commuting variables. Hence the result follows.
	\end{proof}
\end{thm}

\section{Structure of the monoid $M(X,r)$ and its algebra}\label{sec:4}

If $(X,r)$ is a left non-degenerate solution then by Proposition \ref{prop:1} we may (and we shall) identify the
structure monoid $M=M(X,r)$ with its image $\{(a,\phi(a)):a\in A=A(X,r)\}$ in the semidirect product $A\rtimes\mc{G}$,
where $\mc{G}=\mc{G}(X,r)=\gr(\lambda_x\mid x\in X)\s\Sym(X)$, and the map
$\phi\colon A\to\mc{G}$ satisfies $\phi(a)\phi(b)=\phi(a\phi(a)(b))$ for $a,b\in A$.

\begin{lem}\label{lem:4}
	Assume that $(X,r)$ is a finite bijective left non-degenerate solution of the Yang--Baxter equation.
	If $xz=yz$ or $zx=zy$ for some $x,y,z\in M$ then there exists $w\in\Z(M)$ such that $xw=yw$.
	\begin{proof}
		Suppose that $xz=yz$ (the proof in case $zx=zy$ is completely similar) and write $z=(a,\phi(a))$ for some $a\in A$.
		Because $z^n=(a\phi(a)(a)\dotsm\phi(a)^{n-1}(a),\phi(a)^n)$ for each $n\ge 1$, replacing $z$ by some $z^n$ we may assume
		that $\phi(a)=\id$. Moreover, since $a^d\in\Z(A)$ for some $d\ge 1$ (see Lemma \ref{lem:3}) then replacing $z$ by
		$z^d=(a^d,\id)$, we may assume that $a\in\Z(A)$. Since $g(\Z(A))=\Z(A)$ for each $g\in\mc{G}=\mc{G}(X,r)$,
		the element $c=\prod_{g\in\mc{G}}g(a)\in\Z(A)$ is well-defined. It is clear that $g(c)=c$ for each $g\in\mc{G}$.
		Moreover, by induction we prove that $\phi(c^k)=\phi(c)^k$ for $k\ge 1$. Indeed,  
		\[\phi(c^k)=\phi(c^{k-1}c)=\phi(c^{k-1}\phi(c^{k-1})(c))=\phi(c^{k-1})\phi(c)=\phi(c)^{k-1}\phi(c)=\phi(c)^k\]
		for each $k\ge 1$. Hence, replacing $c$ by some $c^k$, we may assume that $\phi(c)=\id$. Define $w=(c,\id)\in M$.
		Clearly $w\in\Z(M)$. Moreover, \[w=(c,\id)=\Big(\prod_{g\in\mc{G}}g(a),\id\Big)=\prod_{g\in\mc{G}}(g(a),\id)=zu,\]
		where $u=\prod_{\id\ne g\in\mc{G}}(g(a),\id)\in A\rtimes\mc{G}$ (note that the element $u$ may not lie in $M$).
		It follows that $xw=yw$, which completes the proof.
	\end{proof}
\end{lem}

As an immediate consequence of Lemma \ref{lem:4} we obtain that the monoid $M$ is left cancellative if and only if it is right cancellative.
Moreover, defining \[\eta_M=\{(x,y)\in M\times M:xz=yz\text{ for some }z\in M\},\] we see that $\eta_M$
is the cancellative congruence of $M$, that is the smallest congruence $\eta$ on $M$ such that the quotient monoid
$M/\eta$ is cancellative. The following proposition gives a description of $\eta_M$ in terms of the cancellative congruence $\eta_A$ of $A=A(X,r)$.

\begin{prop}\label{prop:8}
	Assume that $(X,r)$ is a finite bijective left non-degenerate solution of the Yang--Baxter equation. If $A=A(X,r)$ and $M=M(X,r)$ then
	\[\eta_M=\{((a,\phi(a)),(b,\phi(b))):(a,b)\in\eta_A\text{ and }\phi(a)=\phi(b)\}.\]
	Moreover, there exists $w\in\Z(M)$ and $t\ge 1$ such that \[\eta_M=\{(x,y)\in M\times M:xw^i=yw^i\}\] for all $i\ge t$.
	In particular, the ideal $Mw^t$ is cancellative and if $K$ is a field then $I(\eta_M)=\Ann_{K[M]}(w^i)$ for all $i\ge t$.
	\begin{proof}
		If $(x,y)\in\eta_M$ for some $x=(a,\phi(a))\in M$ and $y=(b,\phi(b))\in M$ then, by the proof of Lemma \ref{lem:4},
		there exists $c\in\Z(A)$ such that $g(c)=c$ for each $g\in\mc{G}=\mc{G}(X,r)$, $\phi(c)=\id$ and $xw=yw$ for $w=(c,\id)\in\Z(M)$.
		Hence \[(ac,\phi(a))=(a,\phi(a))(c,\id)=xw=yw=(b,\phi(b))(c,\id)=(bc,\phi(b)).\] Thus $\phi(a)=\phi(b)$ and $ac=bc$,
		which gives $(a,b)\in\eta_A$. Conversely, if $(a,b)\in\eta_A$ and $\phi(a)=\phi(b)$ then $ac=bc$ for some $c\in\Z(A)$.
		Replacing $c$ by $\prod_{g\in\mc{G}}g(c)\in cA\cap\Z(A)$ we may assume that $g(c)=c$ for each $g\in\mc{G}$. Now,
		\[(a,\phi(a))(c,\phi(c))=(ac,\phi(a)\phi(c))=(bc,\phi(b)\phi(c))=(b,\phi(b))(c,\phi(b)).\]
		Hence $x=(a,\phi(a))\in M$ and $y=(b,\phi(b))\in M$ satisfy $xz=yz$, where $z=(c,\phi(c))\in M$ and thus $(x,y)\in\eta_M$.
		
		To obtain the second equality define $z=\prod_{x\in X}x^d\in\Z(A)$ (here $d\ge 1$ is defined as in Lemma \ref{lem:3}). Let $t\ge 1$
		be such that $\eta_A=\{(a,b)\in A\times A:az^i=bz^i\}$ for all $i\ge t$ (see Proposition \ref{prop:5}).	Since $g(z)=z$ for each
		$g\in\mc{G}$, we get $(z,\phi(z))^n=(z^n,\phi(z)^n)=(z^n,\id)$ for some $n\ge 1$. Define $w=(z^n,\id)\in\Z(M)$. Now, if
		$x=(a,\phi(a))\in M$ and $y=(b,\phi(b))\in M$ then for $i\ge t$ we obtain, by Proposition \ref{prop:5},
		\begin{align*}
			(x,y)\in\eta_M
			& \iff (a,b)\in\eta_A\text{ and }\phi(a)=\phi(b)\\
			& \iff az^{ni}=bz^{ni}\text{ and }\phi(a)=\phi(b)\\
			& \iff xw^i=yw^i,
		\end{align*}
		because
		\begin{equation*}
			xw^i=(a,\phi(a))(z^{ni},\id)=(az^{ni},\phi(a))\qquad\text{and}\qquad yw^i=(b,\phi(b))(z^{ni},\id)=(bz^{ni},\phi(b)).\qedhere
		\end{equation*}
	\end{proof}
\end{prop}

One says that a square-free left non-degenerate solution $(X,r)$ of the Yang--Baxter equation satisfies the so-called
exterior cyclic condition if $r(x,y)=(u,v)$ for some $x,y,u,v\in X$ implies that there exists $z\in X$ such that
$r(v,y)=(u,z)$. This condition was crucial in the study of monoids of I-type (see \cite{JO}). In \cite{Ga2} it is shown
that the exterior cyclic condition holds for a square-free left non-degenerate solution of the Yang--Baxter equation.
Considering the importance of this condition, we include the following generalization of the result in \cite{Ga2},
which can be proved in a similar fashion as in \cite{Ga2}.

\begin{cor}
	Assume that $(X,r)$ is a left non-degenerate solution of the Yang--Baxter equation such that for any $x\in X$
	there exists a unique $y\in X$ satisfying $r(x,y)=(x,y)$. If $x,y,y',u,v,u'\in X$ are such that $r(x,y)=(u,v)$,
	$r(y,y')=(y,y')$ and $r(u,u')=(u,u')$ then there exists $z\in X$ such that $r(v,y')=(u',z)$.
\end{cor}

\begin{lem}\label{lem:new}
    Assume that $(X,r)$ is a left non-degenerate solution of the Yang--Baxter equation. If $\phi$ is the map defined in Proposition \ref{prop:1}
    then the map $\tau\colon X\to X$, defined as $\tau(x)=\phi^{-1}(x)(x)=\lambda_x^{-1}(x)$, satisfies
    \[(x^n,\phi(x^n))=\prod_{i=0}^{n-1}(\tau^i(x),\phi(\tau^i(x)))=(x,\phi(x))(\tau(x),\phi(\tau(x)))\dotsm(\tau^{n-1}(x),\phi(\tau^{n-1}(x)))\]
    for each $x\in X$ and $n\ge 1$. (Here and later by $\tau^0$ we mean the identity map on $X$.)
    \begin{proof}
        We argue by induction on $n$. The formula is clear for $n=1$. Suppose that $n\ge 2$ and assume that the formula holds for $n-1$.
        Then, by the induction hypothesis applied on $\tau(x)$, we get
        \begin{align*}
            (x^n,\phi(x^n))
            & =(x,\phi(x))(\phi^{-1}(x)(x^{n-1}),\phi(\phi^{-1}(x)(x^{n-1})))\\
            & =(x,\phi(x))(\phi(x)^{-1}(x)^{n-1},\phi(\phi(x)^{-1}(x)^{n-1}))\\
            & =(x,\phi(x))(\tau(x)^{n-1},\phi(\tau(x)^{n-1}))\\
            & =(x,\phi(x))(\tau(x),\phi(\tau(x)))\dotsm(\tau^{n-1}(x),\phi(\tau^{n-1}(x))).\qedhere
        \end{align*}
    \end{proof}
\end{lem}

If $(X,r)$ is a left non-degenerate solution of the Yang--Baxter equation then we define the socle
$\Soc(M)$ of the structure monoid $M=M(X,r)$ as \[\Soc(M)=\{(a,\phi(a))\in M:a\in\Z(A)\text{ and }\phi(a)=\id\},\]
where $A=A(X,r)$ is the derived structure monoid of $(X,r)$.

\begin{thm}\label{thm:5}
    Assume that $(X,r)$ is a finite bijective left non-degenerate solution of the Yang--Baxter equation. Let $M=M(X,r)$ and $S=\Soc(M)$.
    Then there exists a finitely generated commutative submonoid $T\s S$ which is normal in $M$ (that is $wT=Tw$ for each $w\in M$),
    and a finite subset $F\s M$ such that $M=\bigcup_{f\in F}Tf$. In particular, if $K$ is a field then $K[M]$ is module-finite normal
    extension of the commutative affine subalgebra $K[T]$. Hence $K[M]$ is a Noetherian PI-algebra. Moreover, if $A=A(X,r)$ then
	\[\clK K[M]=\GK K[M]=\rk M=\clK K[A]=\GK K[A]=\rk A\le|X|\] and the equality holds if and only if the solution $(X,r)$ is involutive.
	\begin{proof}
	    First, we shall show that
	    \begin{equation}
	        \text{there exists }q\ge 1\text{ such that }w^q\in S\text{ for each }w\in M.\label{lab:1}
	    \end{equation}
	    Let $m=|\mc{G}(X,r)|$ and $d\ge 1$ be such that $a^d\in\Z(A)$ for each $a\in A$ (see Lemma \ref{lem:3}).
	    Put $q=md$. If $w=(a,\phi(a))\in M$ then $w^m=(b,\id)$ for some $b\in A$ and thus $w^q=(w^m)^d=(b,\id)^d=(b^d,\id)\in S$.
	    
	    Define $S_0=\{s\in A:(s,\phi(s))\in S\}$. We claim that
	    \begin{equation}
	        \phi(a)(S_0)=S_0\text{ for each }a\in A.\label{lab:2}
	    \end{equation}
	    To prove \eqref{lab:2} choose $s\in S_0$ and define $t=\phi(a)^{-1}(s)\in A$. Since $\phi(a)^{-1}$
	    is an automorphism of $A$ and $s\in\Z(A)$, we get $t\in\Z(A)$. Moreover, $\phi(s)=\id$ leads to
	    \[\phi(a)=\phi(s)\phi(a)=\phi(s\phi(s)(a))=\phi(sa)=\phi(as)=\phi(a\phi(a)(t))=\phi(a)\phi(t),\]
	    which clearly assures that $\phi(t)=\id$. Therefore $\phi(a)^{-1}(s)=t\in S_0$, which means that $\phi(a)^{-1}(S_0)\s S_0$
	    or, equivalently, $S_0\s\phi(a)(S_0)$. Furthermore, because the automorphism $\phi(a)$ is of finite order, we have
	    $\phi(a)=(\phi(a)^{-1})^k$ for some $k\ge 1$. Hence it follows that $\phi(a)(S_0)=(\phi(a)^{-1})^k(S_0)\s S_0$ as well.
	    
	    Next, consider the map $\tau$ defined in Lemma \ref{lem:new}. Since $X$ is a finite set,
	    some power of $\tau$ is an idempotent map, say $\tau^{2p}=\tau^p$ for some $p\ge 1$. We claim that
	    \begin{equation}
	        x^{pq}\in S_0\text{ for each }x\in X.\label{lab:3}
	    \end{equation}
	    To prove \eqref{lab:3} define $v=(\tau^p(x)^p,\phi(\tau^p(x)^p))\in M$. Note that Lemma \ref{lem:new} guarantees that
	    \[v=\prod_{i=0}^{p-1}(\tau^{p+i}(x),\phi(\tau^{p+i}(x))).\]
	    Since $\tau^{2p}=\tau^p$, we get $\tau^{kp}=\tau^p$ for each $k\ge 1$, which yields
	    \[\prod_{i=kp}^{(k+1)p-1}(\tau^i(x),\phi(\tau^i(x)))
	    =\prod_{i=0}^{p-1}(\tau^{kp+i}(x),\phi(\tau^{kp+i}(x)))=\prod_{i=0}^{p-1}(\tau^{p+i}(x),\phi(\tau^{p+i}(x)))=v.\]
	    The above equality together with Lemma \ref{lem:new} and \eqref{lab:1} imply
	    \begin{align*}
	        (\tau^p(x)^{pq},\phi(\tau^p(x)^{pq}))
	        & =\prod_{i=0}^{pq-1}(\tau^i(\tau^p(x)),\phi(\tau^i(\tau^p(x))))=\prod_{i=p}^{p(q+1)-1}(\tau^i(x),\phi(\tau^i(x)))\\
	        & =\Big(\prod_{i=p}^{2p-1}(\tau^i(x),\phi(\tau^i(x)))\Big)\dotsm\Big(\prod_{i=pq}^{p(q+1)-1}(\tau^i(x),\phi(\tau^i(x)))\Big)=v^q\in S.
	    \end{align*}
	    Hence $s=\tau^p(x)^{pq}\in S_0$ and $v^q=(s,\id)$. Let $t=\phi(x^p)^{-1}(x^{pq})\in A$.
	    Because $t=\phi(x^p)^{-1}(x)^{pq}$, we get $t\in\Z(A)$. Moreover, Lemma \ref{lem:new} gives
	    \begin{align*}
	        (x^{p(q+1)},\phi(x^{p(q+1)}))
	        & =\Big(\prod_{i=0}^{p-1}(\tau^i(x),\phi(\tau^i(x)))\Big)\Big(\prod_{j=p}^{p(q+1)-1}(\tau^j(x),\phi(\tau^j(x)))\Big)\\
	        & =(x^p,\phi(x^p))v^q=(x^p,\phi(x^p))(s,\id)=(x^p\phi(x^p)(s),\phi(x^p)),
	    \end{align*}
	    which yields $\phi(x^{p(q+1)})=\phi(x^p)$. Therefore,
	    \begin{align*}
	        \phi(x^p)\phi(t) & =\phi(x^p\phi(x^p)(t))=\phi(x^{p(q+1)})=\phi(x^p)
	    \end{align*}
	    and thus $\phi(t)=\id$. Hence $\phi(x^p)^{-1}(x^{pq})=t\in S_0$, and by \eqref{lab:2} we conclude that
	    $x^{pq}=\phi(x^p)(t)\in S_0$, as claimed.
	    
	    By what we have already shown it follows that $T=\langle(x^{pq},\id)\mid x\in X\rangle$ is a submonoid of the socle $S$.
	    Moreover, as $\phi(a)(x)^{pq}=\phi(a)(x^{pq})\in S_0$ for each $a\in A$ (see \eqref{lab:2} and \eqref{lab:3}), we get
	    \[(a,\phi(a))(x^{pq},\id)=(a\phi(a)(x)^{pq},\phi(a))=(\phi(a)(x)^{pq}a,\phi(a))=(\phi(a)(x)^{pq},\id)(a,\phi(a)).\]
	    Thus, in consequence, $wT=Tw$ for each $w\in M$. Finally, if $X=\{x_1,\dotsc,x_n\}$ with $n=|X|$ then
	    since each element $a\in A$ can be written in the form $a=x_1^{pqk_1+r_1}\dotsm x_n^{pqk_n+r_n}$,
	    where $k_1,\dotsc,k_n\ge 0$ and $0\le r_1,\dotsc,r_n<pq$ (see Remark \ref{rem:1}), we obtain
	    \begin{align*}
	        (a,\phi(a)) & =(x_1^{pq},\id)^{k_1}\dotsm(x_n^{pq},\id)^{k_n}(x_1^{r_1}\dotsm x_n^{r_n},\phi(x_1^{r_1}\dotsm x_n^{r_n})).
	    \end{align*}
	    This clearly leads to a conclusion that $M=\bigcup_{f\in F}Tf$, where
	    \[F=\{(x_1^{r_1}\dotsm x_n^{r_n},\phi(x_1^{r_1}\dotsm x_n^{r_n})):0\le r_1,\dotsc,r_n<pq\}\s M.\]
	    In particular, $K[M]=\sum_{f\in F}K[T]f$. Since $K[T]$ is a commutative affine algebra and the extension $K[T]\s K[M]$
	    is normal (that is $\alpha K[T]=K[T]\alpha$ for each $\alpha\in K[M]$), we get that $K[M]$ is a Noetherian PI-algebra.
	    In this situation, a result of Anan'in \cite{An} implies that $M$ is a linear monoid and then
	    \cite[Proposition 1, p.~221, Proposition 7, p.~280--281, and Theorem 14, p.~284]{O1} yield $\clK K[M]=\GK K[M]=\rk M$.
	    Furthermore, the submonoid $T_0=\{a\in A:(a,\id)\in T\}$ of $A$ is central and it satisfies $x^{pq}\in T_0$ for each
	    $x\in X$ by \eqref{lab:3}. Hence (again referring to Remark \ref{rem:1}) the index of $T_0$ in $A$ is finite (by which we mean that $A$ can be covered by a finite number of cosets of the form $T_0a$
	    for $a\in A$). Therefore, $K[A]$ is a finite module over $K[T_0]$ and thus \[\GK K[M]=\GK K[T]=\GK K[T_0]=\GK K[A].\]
		In view of the last equality, the remaining part of our theorem follows by Theorem \ref{thm:1}.
	\end{proof}
\end{thm}

We finish this section with a positive answer to Conjecture 3.20 posed by Gateva-Ivanova in \cite{Ga2}.

\begin{thm}\label{conditions}
	Assume that $(X,r)$ is a finite bijective left non-degenerate solution of the Yang--Baxter equation.
	Let $M=M(X,r)$. If $K$ is a field then the following conditions are equivalent:
	\begin{enumerate}
		\item $(X,r)$ is an involutive solution.
		\item $M$ is a cancellative monoid.
	    \item $\rk M=|X|$.
		\item $K[M]$ is a prime algebra.
		\item $K[M]$ is a domain.
	    \item $\clK K[M]=|X|$.
	    \item $\GK K[M]=|X|$.
	    \item $\id K[M]=|X|$.
	    \item $K[M]$ has finite global dimension.
	    \item $K[M]$ is an Auslander--Gorenstein algebra.
	    \item $K[M]$ is an Auslander-regular algebra.
	\end{enumerate}
	\begin{proof}
		It is a well-known fact that $(1)\Longrightarrow(9)$ (see \cite[Theorem~1.4]{GVdB}). Moreover,
		$(11)\Longleftrightarrow(9)\Longrightarrow(5)$ and $(10)\Longleftrightarrow(8)\Longrightarrow(7)$ by Theorem \ref{thm:5}
		and \cite[Theorem 1, p.~126]{BG} (see also \cite{SZ}). Clearly $(11)\Longrightarrow(10)$ and $(5)\Longrightarrow(4)$.
		Further, $(4)\Longrightarrow(2)$ follows by Lemma \ref{lem:4}. Since $(1)\Longleftrightarrow(3)\Longleftrightarrow(6)\Longleftrightarrow(7)$
		by Theorem \ref{thm:5}, we obtain the following	diagram of implications
		\[\begin{tikzcd}[row sep=20pt,column sep=20pt]
		    & (3)\ar[d,Leftrightarrow]\ar[r,Leftrightarrow] & (6)\ar[r,Leftrightarrow] & (7) & (8)\ar[l,Rightarrow]\\
		    & (1)\ar[r,Rightarrow] & (9)\ar[r,Leftrightarrow]\ar[d,Rightarrow] & (11)\ar[r,Rightarrow]
		    & (10)\ar[u,Leftrightarrow]\\
		    (2) & (4)\ar[l,Rightarrow] & (5)\ar[l,Rightarrow]
		\end{tikzcd}\]
		Therefore, it is enough to check that $(2)\Longrightarrow(1)$.
		But if $M$ is cancellative and $ca=cb$ for some $a,b,c\in A=A(X,r)$ then
		\begin{align*}
		    (c,\phi(c))(\phi(c)^{-1}(a),\phi(\phi(c)^{-1}(a))) & =(ca,\phi(ca))\\
		    & =(cb,\phi(cb))\\
		    & =(c,\phi(c))(\phi(c)^{-1}(b),\phi(\phi(c)^{-1}(b))).
		\end{align*}
		Hence, by cancellativity of $M$, we get $(\phi(c)^{-1}(a),\phi(\phi(c)^{-1}(a)))=(\phi(c)^{-1}(b),\phi(\phi(c)^{-1}(b)))$.
		Thus $\phi(c)^{-1}(a)=\phi(c)^{-1}(b)$ and $a=b$ follows. Hence $A$ is cancellative and thus $(X,r)$ is an involutive
		solution by Theorem \ref{prop:4}. This finishes the proof.
	\end{proof}
\end{thm}

\begin{rem}
    Moreover, as was proved by Gateva-Ivanova (see \cite[Proposition 3.12]{Ga2}), if the solution $(X,r)$ is square-free then the
    conditions (1)--(11) from Theorem \ref{conditions} are equivalent to Koszulity of the structure algebra $K[M]$.
\end{rem}

Note that in \cite{JVC} it is shown that quadratic monoid is of I-type if and only if it is cancellative and satisfies the cyclic condition.
	   
\section{Prime ideals of $M(X,r)$ and $K[M(X,r)]$}\label{sec:5}

In this section we give a description of certain prime ideals of the algebra $K[M(X,r)]$ over a field $K$ for a
finite bijective left non-degenerate solution $(X,r)$ of the Yang--Baxter equation.
We start with some observations and introduce some notation. As before we make an identification
$M=M(X,r)=\{(a,\phi(a)):a\in A=A(X,r)\}\s A\rtimes\mc{G}$, where $\mc{G}=\mc{G}(X,r)=\gr(\lambda_x\mid x\in X)\s\Sym(X)$,
and the map $\phi\colon A\to\mc{G}$ satisfies $\phi(a)\phi(b)=\phi(a\phi(a)(b))$ for $a,b\in A$.
We first describe all the prime ideals of $M(X,r)$.

Because elements of $A$ are normal, each one-sided ideal of $A$ is a two-sided ideal. For an ideal $I$ of $A$ put \[I^e=\{(a,\phi(a)):a\in I\}.\]
Similarly, if $J$ is an ideal of $M$ then put \[J^c=\{a\in A:(a,\phi(a))\in J\}.\] It is clear that $J^c$ is an ideal of $A$, and $I^e$ is a right
ideal of $M$. Moreover, $I^e$ is an ideal of $M$ if and only if $I$ satisfies $a\phi(a)(I)\s I$ for each $a\in A$ (of course it is enough to consider
$a\in A\setminus I$; let us call such ideals $\phi$-invariant). Thus the rules \[I\mapsto I^e\qquad\text{and}\qquad J\mapsto J^c\]
define mutually inverse bijections (actually mutually inverse lattice isomorphisms) between the set consisting of all $\phi$-invariant ideals
of $A$ and the set consisting of all ideals of $M$.

\begin{lem}\label{lem:5}
	Assume that $(X,r)$ is a finite left non-degenerate solution of the Yang--Baxter equation.
	If $P$ is a prime ideal of $M=M(X,r)$ then $P=I^e$ with $I$ a semiprime ideal of $A=A(X,r)$.
	Thus $P=(Q_1\cap\dotsb\cap Q_r)^e$ for some prime ideals $Q_1,\dotsc,Q_r$ of $A$ that are minimal over $I$.
	\begin{proof}
		Let $P=I^e$ be a prime ideal of $M$. We need to prove that $I$ is a semiprime ideal of $A$. To do so, assume $J$ is an ideal of $A$ that
		contains $I$ and such that $J/I$ is nil. We claim that $I=J$. First we show that the right ideal $J^e=\{(j,\phi(j)):j\in J\}$ of $M$ is nil
		modulo $P$. Indeed, take $x=(j,\phi(j))\in J^e$. For any $n\ge 1$ we have that \[x^n=(j\phi(j)(j)\dotsm\phi(j)^{n-1}(j),\phi(j)^n).\]
		Hence for a large enough $n\ge 1$, we have $x^n=(y,\id)$ with $y\in J$. And thus for some $m\ge 1$, we get that $x^{nm}=(y^m,\id)\in I^e$
		and $y^m\in I$. This proves that $J^e$ is indeed nil modulo $P$. Hence, $J^e/P$ is nil submonoid of the monoid $M/P$.
		Since $M$ and thus also $M/P$ satisfies the ascending chain condition, it is well-known
		(cf. \cite[Proposition 17.22]{Fa} or \cite[Theorem 2.4.10]{JO}) that $J^e/P$ is nilpotent.
		Since $P$ is a prime ideal we get that $J^e\s P=I^e$ and thus $J=I$, as desired.
	\end{proof}
\end{lem}

With notation as above, since $P$ is a left ideal we have that $a\phi(a)(Q_1\cap\dotsb\cap Q_r)\s Q_1\cap\dotsb\cap Q_r$ for every $a\in A$.
As $\phi(a)\in\Aut(A)$ this condition is equivalent with
\begin{equation}\label{pr}
	\text{for every }1\le i \le r\text{ and for every }a\in A\setminus Q_i\text{ there exists }1\le j\le r\text{ such that }\phi(a)(Q_j)\s Q_i.
\end{equation}
Renumbering, if necessary, we may assume that $Q_1,\dotsc,Q_k$ are all the prime ideals of least height among all primes $Q_1,\dotsc,Q_r$.
Then, condition \eqref{pr} yields that for every $1\le i\le k$ and for every $a\in A\setminus Q_i$ there exists $1\le j\le k$ such that
$\phi(a)(Q_j)=Q_i$. Hence $(Q_1\cap\dotsb\cap Q_k)^e$ is an ideal of $M$. We claim that $k=r$. Suppose the contrary, i.e., suppose $k<r$.
First note that $(Q_{k+1}\cap\dotsb\cap Q_r)^e$ and $(Q_1\cap\dotsb\cap Q_k)^e$ are right ideals of $M$.
Furthermore,
 \begin{align*}
	(Q_{k+1}\cap\dotsb\cap Q_r)^e(Q_1\cap\dotsb\cap Q_k)^e &\s\bigcup_{a\in Q_{k+1}\cap\dotsb\cap Q_r}
	(a\phi(a)(Q_1\cap\dotsb\cap Q_k),\phi(a\phi(a)(Q_1\cap\dotsb\cap Q_k)))\\
	&\s(Q_1\cap\dotsb\cap Q_k\cap Q_{k+1}\cap\dotsb\cap Q_r)^e=P.
 \end{align*}
Thus, $(Q_{k+1}\cap\dotsb\cap Q_r)^e\s P$ or $(Q_1\cap\dotsb\cap Q_k)^e\s P$. The former would imply
that $Q_s\s Q_1$ for some $k+1\le s\le r$. Hence, since all the primes involved are minimal over $I$,
we would get $Q_s=Q_1$, a contradiction.

Hence, we have proved the first part of the following lemma. The second part is then a translation of condition \eqref{pr}.

\begin{lem}\label{lem:6}
	Assume that $(X,r)$ is a finite left non-degenerate solution of the Yang--Baxter equation.
	If $P$ is a prime ideal of $M=M(X,r)$ then $P=I^e$ with $I$ a semiprime ideal of $A=A(X,r)$.
	Thus $P=(Q_1\cap\dotsb\cap Q_r)^e$ where $Q_1,\dotsc,Q_r$ are prime ideals of $A$ all of the same height, and furthermore,
	\[\text{for every }1\le i \le r\text{ and for every }a\in A\setminus Q_i\text{ there exists }1\le j\le r\text{ such that }\phi(a)^{-1}(Q_i)=Q_j.\]
	The set of prime ideals $\{Q_1,\dotsc,Q_r\}$ will be denoted as $\Spec(P)$. Consequently,
	\[\text{if }Q\in\Spec(P)\text{ then }\{\phi(a)^{-1}(Q):a\in A\setminus Q\}\s\Spec(P).\]
\end{lem}

We now focus on the converse process and investigate whether for a prime ideal $Q$ of $A$ there exists a prime ideal $P$ of $M$
such that $Q\in\Spec(P)$. To do so we recursively introduce some sets $\mc{S}_n=\mc{S}_n(Q)$ consisting of prime ideals of $A$.
Put \[\mc{S}_1=\mc{S}_1(Q)=\{Q\}\] and \[\mc{S}_{n+1}=\mc{S}_{n+1}(Q)=\{\phi(a)^{-1}(Q'):Q'\in\mc{S}_n\text{ and }a\in A\setminus Q'\}.\]
Since $1\in A\setminus Q'$ for $Q'\in\mc{S}_n$ and $\phi(1)=\id$ we get $\mc{S}_n\s\mc{S}_{n+1}$. Because $A$ has only finitely many
prime ideals there exists $n=n(Q)\ge 1$ such that $\mc{S}_i=\mc{S}_n$ for all $i\ge n$. Put \[P(Q)=\bigcap_{Q'\in\mc{S}_n}Q'.\]
 
We also need the following lemma.

\begin{lem}\label{lem:7}
	Assume that $(X,r)$ is a finite bijective left non-degenerate solution of the Yang--Baxter equation.
	Then there exists $t\ge 1$ such that $\phi(a^t)=\id$ for each $a\in A=A(X,r)$.
	\begin{proof}
		From the proof of Theorem \ref{thm:5} we know that there exists $d\ge 1$ such that $a^d\in\Z(A)$
		for each $a\in A$ and $\phi(x^d)=\id$ for each $x\in X$ (it is enough to take $d=pq$, where $p,q$ are as in the proof of
		Theorem \ref{thm:5}; note also that $d$ used here is a multiple of $d$ from Lemma \ref{lem:3}).
		Let $C=\free{x^d\mid x\in X}$, a submonoid of $\Z(A)$.
		Let $\sim$ denote the equivalence relation on $A$ defined by $a_1\sim a_2$ if $c_1a_1=c_2a_2$ for some $c_1,c_2\in C$.
		Because $C$ is a central submonoid of $A$ we have that $\sim$ is a congruence on $A$. Denote by $\ov{a}$ the natural
		image of $a\in A$ in the monoid $\ov{A}=A/\!\!\sim$. Clearly $\ov{A}=\free{\ov{x}\mid x\in X}$.
		As $\ov{x}^d=\ov{x^d}=\ov{1}$, the monoid $\ov{A}$ is a group and, by Remark \ref{rem:1},
		it follows that $\ov{A}$ is a finite group, say of order $t\ge 1$. Then, for every $a\in A$, we obtain that $\ov{a}^t=\ov{1}$.
		Hence, for every $a\in A$ there exist $c_1,c_2\in C$ such that $c_1a^t=c_2$. Since $\phi(c_1)=\phi(c_2)=\id$, we conclude that
		\begin{equation*}
			\phi(a^t)=\phi(c_1)\phi(a^t)=\phi(c_1\phi(c_1)(a^t))=\phi(c_1a^t)=\phi(c_2)=\id,
		\end{equation*}
		as desired.
	\end{proof}
\end{lem}

\begin{lem}
	Assume that $(X,r)$ is a finite bijective left non-degenerate solution of the Yang--Baxter equation.
	If $Q$ is a prime ideal of $A=A(X,r)$ then $P(Q)^e$ is a prime ideal of $M=M(X,r)$.
	\begin{proof}
		First we show that $P(Q)^e$ is an ideal of $M$. For this we need to show that condition \eqref{pr} holds for the set of primes $\mc{S}_n$
		(where $n=n(Q)$). So, let $Q'\in\mc{S}_n$ and $a\in A\setminus Q'$. Then by the definition of $\mc{S}_n$ we have that
		$\phi(a)^{-1}(Q')\in\mc{S}_{n+1}=\mc{S}_n$ and thus condition \eqref{pr} follows.
		
		Second we prove that $P(Q)^e$ is a prime ideal of $M$. To do so, consider
		\[\mc{F}=\{I^e:I\text{ is an ideal of }A\text{ such that }I\s Q\text{ and }I^e\text{ is an ideal of }M\}.\]
		By the first part $P(Q)^e\in\mc{F}$ and thus $\mc{F}\ne\vn$. Because of Zorn's Lemma, there exists a maximal
		(for the inclusion relation) element of $\mc{F}$, say $I^e$. We claim that $I^e$ is prime ideal of $M$.
		To prove this, suppose $J^e$ and $K^e$ are ideals of $M$, with $J$ and $K$ ideals of $A$ that properly contain $I$, such that $J^eK^e\s I^e$.
		Then, because of the maximality, there exist $j\in J\setminus Q$ and $k\in K\setminus Q$ and $(j,\phi(j))M(k,\phi(k))\s I^e$.
		Because of Lemma~\ref{lem:7}, let $d\ge 1$ be such that $\phi(j^d)=\id$. Then,
		\[(j,\phi(j))\cdot(\phi(j)^{-1}(j^{d-1}),\phi(\phi(j)^{-1}(j^{d-1})))\cdot(k,\phi(k))=
		(j^d,\id)\cdot(k,\phi(k))=(j^dk,\phi(k))\in I^e.\]
		Hence $j^dk\in I\s Q$, in contradiction with $Q$ being a prime ideal in the monoid $A$ that consists of normal elements.
		So, indeed $I^e$ is a prime ideal of $M$.
		
		Hence, by Lemmas \ref{lem:5} and \ref{lem:6} we know that $I=Q_1\cap\dotsb\cap Q_r$, an intersection
		of primes ideals of $A$ of the same height, and $\mc{S}_n(Q)\s\Spec(I^e)$. So $I\s P(Q)\s Q$ and thus $I^e\s P(Q)^e$.
		Since $P(Q)^e\in\mc{F}$, the maximality condition yields that $I^e=P(Q)^e$ and the result follows.
	\end{proof}
\end{lem}

The previous lemmas together with results from Section \ref{sec:3} give a full description of the prime ideals in $M(X,r)$.

\begin{prop}
	Assume that $(X,r)$ is a finite bijective left non-degenerate solution of the Yang--Baxter equation.
	The prime ideals of $M=M(X,r)$ are precisely the ideals $P(Q)^e$, where $Q$ runs through the prime ideals of $A=A(X,r)$.
	Further, \[P(Q)=\bigcap_{Q'\in\mc{S}_{n(Q)}(Q)}Q'\] is an intersection of prime ideals of $A$ of the same height and
	$\h P(Q)^e=\h Q$. In particular, the map \[\Spec(A)\to\Spec(M)\colon Q\mapsto P(Q)^e\] satisfies going-up, going-down and incomparability.
	\begin{proof}
		The first part has been proved. The second part follows now at once.
	\end{proof}
\end{prop}

We also have the following analog of the second part of Proposition \ref{prop:7}. That is, prime ideals of the algebra
$K[M(X,r)]$ over a field $K$ not intersecting the monoid $M(X,r)$ are determined by prime ideals of the group algebra $K[G(X,r)]$.

\begin{prop}\label{prop:11}
	Assume that $(X,r)$ is a finite bijective left non-degenerate solution of the Yang--Baxter equation.
	Let $M=M(X,r)$. If $K$ is a field then there exists an inclusion preserving bijection between the
	set of prime ideals $P$ of $K[M]$ satisfying $P\cap M=\vn$ and the set of all prime ideals of the
	group algebra $K[G]$, where \[G=G(X,r)=\gr(X\mid xy=\lambda_x(y)\rho_y(x)\text{ for all }x,y\in X).\]
	Moreover, the cancellative monoid $\ov{M}=M/\eta_M$ has a group of quotients, which is  equal to
	the central localization $\ov{M}\free{z}^{-1}$ for some $z\in\Z(\ov{M})$, and $G\cong\ov{M}\free{z}^{-1}$.
	\begin{proof}
		By Lemma \ref{lem:3} there exists $d\ge 1$ such that $\phi(a)^d=\id$ and $a^d\in\Z(A)$
		for each $a\in A=A(X,r)$. Moreover, if $x\in X$ then $(x,\phi(x))^{d^2}=(y,\id)^d=(a_x,\id)$, where
		$y=x\phi(x)(x)\dotsm\phi(x)^{d-1}(x)\in A$ and $a_x=y^d\in\Z(A)$. Define
		$c_x=\prod_{g\in\mc{G}}g(a_x)$, where $\mc{G}=\mc{G}(X,r)$. Clearly $c_x\in\Z(A)$ and $g(c_x)=c_x$
		for each $g\in\mc{G}$. If $b_x=\prod_{\id\ne g\in\mc{G}}g(a_x)\in\Z(A)$ then $a_xb_x=c_x$ and thus
		\[(x,\phi(x))^{d^2}(b_x,\phi(b_x))=(a_x,\id)(b_x,\phi(b_x))=(a_xb_x,\phi(b_x))=(c_x,\phi(c_x)).\]
		Moreover, $\phi(c_x^d)=\phi(c_x)^d=\id$. Hence $[(x,\phi(x))^{d^2}(b_x,\phi(b_x))]^d=(c_x,\phi(c_x))^d=(c_x^d,\id)$,
		and it follows that if \[z=\prod_{x\in X}(c_x^d,\id)=(\prod_{x\in X}c_x^d,\id)\] then $z\in\Z(M)$ and
		the central localization $\ov{M}\free{z}^{-1}$ (here by $z$ we understand the image of $z\in M$ in $\ov{M}$)
		is equal to the group of quotients of $\ov{M}$. Finally, the remaining part of the proof is similar to the proof
		of Proposition \ref{prop:7}. Thus the result follows.
	\end{proof}
\end{prop}

\section{Divisibility in $M(X,r)$}\label{sec:6}

In the previous section we have shown that prime ideals are sets that are determined by divisibility of some generators.
In this section, we go deeper into this. This has been done earlier, but for different quadratic monoids, in several
papers to prove that the algebra is Noetherian and PI, which we already know.  As before, throughout this section $(X,r)$
denotes a finite bijective left non-degenerate solution of the Yang--Baxter equation. Let us now relate the structure of
$M=M(X,r)$ to substructures determined by divisibility by generators. Since each element of
$A=A(X,r)$ is a normal element, left divisibility in $A$ by an element is the same as right divisibility by that element.
In the monoid $M$ we will have to use the terminology left and right divisible. Let $|X|=n$. For $1\le i\le n$ put
\[A_i=\{a\in A:a\text{ is divisible by at least }i\text{ elements of }X\}\] and
\[M_i=\{m\in M:m\text{ is left divisible by at least }i\text{ elements }(x,\phi(x))\text{ with }x\in X\}.\]
Clearly, \[M_i=A_i^e=\{(a,\phi(a)):a\in A_i\}.\] Since $A_i$ is a $\phi$-invariant ideal of $A$, it follows that each $M_i$
is an ideal of $M$. Hence we get an ideal chain in $M$: \[M_n\s M_{n-1}\s\dotsb\s M_1\s M.\]
Note that the equality $M_i=M_{i+1}$ is possible, for instance the structure monoid $M=M(X,r)$ of the solution $(X,r)$
defined in Example~\ref{ex:1} satisfies $M_1=M_2$.

The following lemmas, propositions and proofs are completely analogous to those for monoids of skew type given in \cite{JO}.
We have included them for completeness' sake.

For a non-empty subset $Y\s X$ define
\[M_Y=\bigcap_{y\in Y}(y,\phi(y))M\qquad\text{and}\qquad D_Y=M_Y\setminus\bigcup_{x\in X\setminus Y}M_{\{x\}}.\]
Clearly, the set $M_Y$ consists of all elements of $M$ that are left divisible by generators $(y,\phi(y))$ with $y\in Y$,
and the set $D_Y$ consists of all elements of $M$ that are precisely left divisible by those generators.
Obviously, for each $1\le i\le n$ we have \[M_i=\bigcup_{Y\s X,\,|Y|=i}M_Y.\]

The following lemma is clear, by using the fact that $M_i=A_i^e$ and $A_i$ are $\phi$-invariant ideals of $A$.

\begin{lem}[{cf. \cite[Theorem 9.3.7]{JO}}]
	Assume that $(X,r)$ is a finite left non-degenerate solution of the Yang--Baxter equation. If $M=M(X,r)$ and $n=|X|$ then
	\[M_X=M_n\s M_{n-1}\s\dotsb\s M_1\s M\] is a chain of ideals in $M$. 
\end{lem}

The following technical lemma will prove to be crucial in the proof of the main result of this section.
It proves that under certain conditions we can show left divisibility by words.

\begin{lem}[{cf. \cite[Lemma 9.3.8]{JO}}]\label{lem:10}
	Assume that $(X,r)$ is a finite bijective left non-degenerate solution of the Yang--Baxter equation.
	Let $M=M(X,r)$ and $n=|X|$. Suppose that $Y\s X$ and $|Y|=i$, where $1\le i\le n$. If $\vn\ne Z\s Y$
	and $s\in D_Z$ satisfies $|s|=k$ then \[M_i^k\cap D_Y\s sM.\]
	\begin{proof}
		If $k=1$ then the claim is obvious. So assume $k\ge 2$. To shorten the notation put \[s_x=(x,\phi(x))\in M\] for $x\in X$
		and write $s=s_{x_1}\dotsm s_{x_k}$ with $x_1,\dotsc,x_k\in X$. Let $a=a_1\dotsm a_k\in D_Y$,
		where $a_1,\dotsc,a_k\in M_i$. Since $D_Y\s M_i\setminus M_{i+1}$ and $M_{i+1}$ is an ideal of $M$, if non-empty,
		it is clear that $a_1,\dotsc,a_k\in M_i\setminus M_{i+1}$.		
		As $s\in D_Z$ and $Z\s Y$, it follows that $x_1\in Y$. Because $a\in D_Y$, we obtain $a_1\in D_Y$. Hence, there exists $b_1\in M$ 
		such that $a_1=s_{x_1}b_1$. Thus, $a_1a_2=s_{x_1}c_1$, where $c_1=b_1a_2$. Clearly, $a_1a_2\in M_i\setminus M_{i+1}$, which
		implies that $c_1\in M_i\setminus M_{i+1}$. Suppose we have shown that
		\[a_1\dotsm a_j=s_{x_1}\dotsm s_{x_{j-1}}c_{j-1}\] for some $1<j<k$ and $c_{j-1}\in M_i\setminus M_{i+1}$.
		We claim that $c_{j-1}\in s_{x_j}M$. Let $W\s X$ be such that $|W|=i$ and $c_{j-1}\in D_W$. Consider the set
		\[U=\{x\in X:s_{x_1}\dotsm s_{x_{j-1}}s_x\in D_V\text{ for some }V\s Y\}.\] As $(X,r)$ is left non-degenerate, it follows that $|U|\le |Y|=i$.
		Since $a\in D_Y$, it follows that $a_1\dotsm a_j\in D_Y$. Because $c_{j-1}\in D_W$, we obtain that $W\s U$. Thus $|U|=i$ and $W=U$.
		Since $s_{x_1}\dotsm s_{x_j}$ is a left initial segment of $s\in D_Z$ and $Z\s Y$, we also get that $x_j\in U=W$.
		As $c_{j-1}\in D_W$, it follows that $c_{j-1}\in s_{x_j}M$, as claimed.
		
		Now, write $c_{j-1}=s_{x_j}b_j$ for some $b_j\in M$. Then \[a_1\dotsm a_ja_{j+1}=s_{x_1}\dotsm s_{x_{j-1}}s_{x_j}b_ja_{j+1}.\]
		Define $c_j=b_ja_{j+1}$. Then $c_j\in M_i\setminus M_{i+1}$. Thus, by induction, we obtain that
		$a=a_1\dotsm a_k\in s_{x_1}\dotsm s_{x_k}M=sM$, and the result is shown.
	\end{proof}
\end{lem}

Recall that by \[I(\eta)=\Span_K\{x-y:(x,y)\in \eta\},\] the $K$-linear span of the set consisting of all elements $x-y$ with $(x,y)\in\eta$,
we understand the ideal of the algebra $K[M]$ associated to a congruence $\eta$ on the monoid $M$. Moreover, $\eta_M$ denotes the cancellative
congruence of $M$ (see Proposition \ref{prop:8} and the comment above).

\begin{lem}\label{lem:11}
	Assume that $(X,r)$ is a finite bijective left non-degenerate solution of the Yang--Baxter equation.
	Let $M=M(X,r)$. If $K$ is a field then \[I(\eta_M)=\Ann_{K[M]}(w^m)=\Ann_{K[M]}(M_X^m)\] for some $w\in M_X\cap\Z(M)$ and some $m\ge 1$.
	\begin{proof}
		Let $A=A(X,r)$. Define $z=\prod_{x\in X}x^d\in\Z(A)$ and $w=(z^n,\id)\in\Z(M)$ as in the proof of Proposition \ref{prop:8}.
		Clearly $w\in M_X$. Therefore, Lemma \ref{lem:10} implies that $M_X^k\s wM\s M_X$, where $k=|w|$. Now, if $t\ge 1$ is such that
		$I(\eta_M)=\Ann_{K[M]}(w^i)$ for all $i\ge t$ (see Proposition \ref{prop:8}) then, since $w^{kt}M\s M_X^{kt}\s w^tM$, we conclude that
		\[I(\eta_M)=\Ann_{K[M]}(w^t)\s\Ann_{K[M]}(M_X^{kt})\s\Ann_{K[M]}(w^{kt})=I(\eta_M).\]
		This shows the result with $m=kt$.
	\end{proof}
\end{lem}

The following proposition provides us information on prime ideals $P$ of the algebra $K[M(X,r)]$, which intersect the monoid $M(X,r)$
non-trivially.

\begin{prop}[{cf. \cite[Proposition 9.5.3]{JO}}]\label{prop:12}
	Assume that $(X,r)$ is a finite bijective left non-degenerate solution of the Yang--Baxter equation. Let $M=M(X,r)$. If $K$ is a field
	and $P$ is a prime ideal of $K[M]$ such that $P\cap M\ne\vn$ then \[P\cap M=\bigcup_{Y\in\mc{F}}D_Y,\]
	where $\mc{F}=\{Y\s X:Y\ne\vn\text{ and }D_Y\cap P\ne\vn\}$. Moreover, if $Y\in\mc{F}$ and $Y\s Z\s X$ then $D_Z\s P$.
	\begin{proof}
		The inclusion $P\cap M\s\bigcup_{Y\in\mc{F}}D_Y$ is obvious. We prove the reverse inclusion by contradiction.
		So, suppose that there exists $Y\in\mc{F}$ such that $D_Y\nsubseteq P$. Choose such a set $Y$ with maximal $i=|Y|$.
		We claim that $i<|X|$. Indeed, if $s\in P\cap M$ then $ws\in P\cap M_X$ (we use the notation introduced in the
		proof of Lemma \ref{lem:11}). Therefore, by Lemma \ref{lem:10}, we get $M_X^n\s wsM\s P$, where $n=|ws|$.
		Hence $D_X=M_X\s P$ and thus $Y\ne X$, as claimed. Let $a\in D_Y\cap P$ and set $k=|a|$. Consider an arbitrary subset $Z\s X$
		such that $Y\s Z$ and $D_Z\ne\vn$. By Lemma \ref{lem:10},
		\begin{equation}
			D_Z^k\cap D_Z\s aM\s P.\label{DZ}
		\end{equation}
		If $Z\ne Y$ and $D_Z^k\cap D_Z\ne\vn$, then $|Z|>i$ and $D_Z\cap P\ne\vn$ by \eqref{DZ}. Hence, by the definition of $i$,
		we get $D_Z\s P$. Whereas, if $Z\ne Y$ and $D_Z^k\cap D_Z=\vn$, then $D_Z^k\s\bigcup_{Z\sn V}D_V$. In the latter case
		the given argument can be applied to every $V\s X$ such that $Z\sn V$ and $D_V\ne\vn$. Continuing this process, after
		a finite number of steps, we obtain that $I=\bigcup_{Y\sn Z}D_Z$ is nilpotent modulo $P$. Since $I$ is a right ideal of
		$M$ and because $P\cap M$ is a prime ideal of $M$, it follows that $I\s P$. Applying \eqref{DZ} to $Z=Y$,
		we thus have proved that \[D_Y^k\s(D_Y^k\cap D_Y)\cup I\s P.\]
		As $D_Y\cup I=\bigcup_{Y\s Z}D_Z$ is a right ideal of $M$ and it is nilpotent modulo $P$, we conclude that $D_Y\s P$,
		a contradiction. The second part of the result follows from the proof above.
	\end{proof}
\end{prop}

We now are in a position to prove the main result of this section.

\begin{thm}[{cf. \cite[Proposition 9.5.2]{JO}}]
	Assume that $(X,r)$ is a finite bijective left non-degenerate solution of the Yang--Baxter equation. If $M=M(X,r)$ and $K$ is a field
	then the following properties hold:
	\begin{enumerate}
		\item $I(\eta_M)=\Ann_{K[M]}(w^m)=\Ann_{K[M]}(M_X^m)$ for some $w\in M_X\cap\Z(M)$ and some $m\ge 1$.
		\item $I(\eta_M)\s P$ for any prime ideal $P$ of $K[M]$ such that $P\cap M=\vn$.
		\item $M_X\s P$ for any prime ideal $P$ of $K[M]$ such that $P\cap M\ne\vn$. In particular, $w\in P$.
		\item There exists at least one minimal prime ideal $P$ of $K[M]$ such that $P\cap M=\vn$.
		\item If $\ch K=0$ then
		\[\mc{J}(K[M])=\mc{B}(K[M])=I(\eta_M)\cap\bigcap_{P\in\mc{P}}P,\]
		where $\mc{P}=\{P\in\Spec(K[M]):P\cap M\ne\vn\}=\{P\in\Spec(K[M]):w\in P\}$.
	\end{enumerate}
	\begin{proof}
		(1), (2) and (3) follow from Lemma \ref{lem:11}, Proposition \ref{prop:8} and Proposition \ref{prop:12}, respectively. 
		Clearly $M\cap\mc{B}(K[M])=\vn$ (elements of $M$ are not nilpotent). Therefore, $w\notin\mc{B}(K[M])$ and (4) also follows.
		If $\ch K=0$ then the algebra $K[M/\eta_M]\cong K[M]/I(\eta_M)$ is semiprime (see \cite[Theorem 3.2.8]{JO}).
		Hence $\mc{B}(K[M])\s I(\eta_M)$ and $I(\eta_M)$ is equal to the intersection of all prime ideals $P$ of $K[M]$ such
		that $I(\eta_M)\s P$. Thus the second equality in (5) follows. Since the Jacobson radical of an affine
		PI-algebra is nilpotent (cf. \cite{Br}), it equals the prime radical, which ends the proof.
	\end{proof}
\end{thm}

\section{Prime images of $K[M(X,r)]$.}\label{sec:7}

Assume that $(X,r)$ is a finite bijective left non-degenerate solution of the Yang--Baxter equation. Let $M=M(X,r)$ and let $K$ be a field.
The aim in this section is to provide a matrix-type representation of the prime algebra $K[M]/P$ for each prime ideal $P$ of $K[M]$. 
We do this by showing that the classical ring of quotients $\Qcl(K[M]/P)$ is the same as $\Qcl(\M_v(K[G]/P_0))$, where $P_0$ is
a prime ideal of a group algebra $K[G]$ with $G$ the group of quotients of a cancellative subsemigroup of $M$ and $v\ge 1$ is determined
by the number of orthogonal cancellative subsemigroups of an ideal in $M/(P\cap M)$. If $P$ is such that $P\cap M=\vn$ then this has been shown
in Proposition \ref{prop:11}. Hence, in the remainder of this section we assume that $P$ is a prime ideal of $K[M]$ with $P\cap M\ne\vn$.
Note that $K[M]/P$ is an epimorphic image of the contracted monoid algebra $K_0[M/(P\cap M)]$. Hence we determine a representation of
\[S=M/(P\cap M).\] 

As a first step we make use of a result of Anan'in (see \cite[Theorem 3.5.2]{JO} or \cite{An}) that yields that the Noetherian
PI-algebra $K_0[S]$ embeds into a matrix algebra $\M_m(L)$ over a field extension $L$ of $K$. Thus, we will consider $S$
as a submonoid of the multiplicative monoid $\M_m(L)$. By \cite[Proposition 5.1.1]{JO} (and the fact that $S$ satisfies the ascending
chain condition on left and right ideals) it follows that $S$ intersects non-trivially finitely many $\mc{H}$-classes of $\M_m(L)$
(i.e., the maximal subgroups of $\M_m(L)$), say $G_1,\dotsc,G_k$. Since $K_0[S]$ is a PI-algebra, also each $K[S\cap G_i]$ is a PI-algebra.
Hence, $S\cap G_i$ has a group of quotients $\gr(S\cap G_i)$ which is abelian-by-finite (cf. \cite[Theorem 3.1.9]{JO}).

For every $1\le i\le k$, let $e_i$ denote the idempotent of the maximal subgroup $G_i$ and fix $s_i\in S\cap G_i$. Because of Lemma \ref{lem:3}
we may choose $s_i$ in the center of $A(X,r)$. So, $e_i=s_is_i^{-1}$, where $s_i^{-1}$ denotes the inverse of $s_i$ in $G_i$. Without loss
of generality, we may assume that $s_i=(a_i,\id) $ with $a_i$ in the center of $A(X,r)$. One then proves as in \cite[Lemma 2.4]{JVC}
that $e_ie_j=e_je_i$ for all $1\le i,j\le k$. Hence, \[\free{e_1,\dotsc,e_k}\cup\{\theta\}=\{e_1,\dotsc,e_k\}\cup\{\theta\}\]
is an abelian semigroup (where $\theta$ is the zero element of $S$). By \cite[Theorem 3.5]{O2} we get that the linear semigroup
$S$ has an ideal chain 
\begin{equation}\label{IC}
	S_0\s T_1\s S_1\s\dotsb\s S_{m-1}\s T_m\s S_m=S
\end{equation}
with each \[N_j=T_j/S_{j-1}=(T_j\setminus S_{j-1})\cup\{\theta\}\text{ a nilpotent ideal of }S/S_{j-1}\]
(and it actually is a union of nilpotent ideals of nilpotency index $2$) and each
\[S_j/T_j=(S_j\setminus T_j)\cup\{\theta\}\s\mc{M}_j/\mc{M}_{j-1}\] a $0$-disjoint union of uniform subsemigroups (for the terminology see
\cite[Section 2.2]{JO}), say $U_{\alpha}^{(j)}$ (with $\alpha$ in an indexing set $\mc{A}_j$), of $\mc{M}_j/\mc{M}_{j-1}$ that intersect
different $\mc{R}$-classes and different $\mc{L}$-classes of $\mc{M}_j/\mc{M}_{j-1}$ (here $\mc{M}_j$ denotes the ideal in $\M_m(L)$ consisting
of matrices of rank at most $j$). Recall that it is well-known that $\mc{M}_j$ are the only ideals of the multiplicative monoid $\M_m(L)$ and 
$\mc{M}_j/\mc{M}_{j-1}$ is a completely $0$-simple semigroup with maximal subgroups isomorphic to $\GL_j(L)$. Moreover,
\[\text{each } N_j\text{ does not intersect }\mc{H}\text{-classes of }\mc{M}_j/\mc{M}_{j-1}\text{ intersected by }S_j\setminus T_j\]
and \[U_{\alpha}^{(j)}U_{\beta}^{(j)}\s N_j\text{ for all }\alpha\ne\beta
\text{ and }U_{\alpha}^{(j)}N_jU_{\alpha}^{(j)}=\{\theta\}\text{ in }\mc{M}_j/\mc{M}_{j-1}.\]
In particular, each $U_{\alpha}^{(j)}$ can be considered as an ideal in $S/T_j$.

Because $S$ is a prime monoid with $0$-element, it follows that the lowest non-zero ideal in the chain \eqref{IC}
is of the type $S_j$ (i.e., $T_j=\{\theta\}$). So $N_j=\{ \theta\}$ and thus $S_j\s\mc{M}_j/\mc{M}_{j-1}$ is a $0$-disjoint
union of the uniform subsemigroups $U_{\alpha}^{(j)}$, and each $U_{\alpha}^{(j)}$ is an ideal of $S$.
As $U_{\alpha}^{(j)} U_{\beta}^{(j)}\s N_j=\{\theta\}$ for $\alpha\ne\beta$, and because $S$ is a prime monoid we get that
$S_j=U_{\alpha}^{(j)}$ for some $\alpha$, and it is a uniform subsemigroup of the completely $0$-simple semigroup $\mc{M}_j/\mc{M}_{j-1}$.
Renumbering $G_1,\dotsc,G_k$, if necessary, we may assume that $G_1,\dotsc,G_v$ are all the maximal subgroups
of $\mc{M}_j$ that intersect $S$ non-trivially. So, for each $1\le r\le v$, the semigroup $S\cap G_r$ is cancellative.

We also know that $S_j=U_{\alpha}^{(j)}$ is contained in the smallest completely $0$-simple subsemigroup $\wh{U}_\alpha^{(j)}$ of 
$\mc{M}_j/\mc{M}_{j-1}$ (see, e.g., \cite[Proposition 2.2.1]{JO}). That is, $U_{\alpha}^{(j)}$ intersects all nonzero $\mc{H}$-classes
of $\wh{U}_\alpha^{(j)}$ and every maximal subgroup $H$ of $\wh{U}_\alpha^{(j)}$ is generated by $U_{\alpha}^{(j)}\cap H$
(so $H=\gr(S\cap G_i)$ for some $i$ and $\gr(S\cap G_1)\cong\dotsb\cong\gr(S\cap G_v)$ is an abelian-by-finite group).
 
To simplify notation, we write $U_{\alpha}^{(j)}$ as $U$ and $\wh{U}_\alpha^{(j)}$ as $\wh{U}$. By the above, the idempotents of $\wh{U}$ commute.
Since $\wh{U}$ is completely $0$-simple, this implies that these idempotents are pairwise orthogonal. Since $S$ intersects non-trivially only
finitely $\mc{H}$-classes of $\mc{M}_j/\mc{M}_{j-1}$, the completely $0$-simple semigroup $\wh{U}$ has only finitely many rows and columns.
It follows that the sandwich matrix of $\wh{U}$ contains precisely one non-zero element in each row and column. So, reindexing if necessary,
we may assume that the sandwich matrix is a diagonal matrix, and thus also $\wh{U}$ has the same number of rows and columns. It is then well-known
(see, e.g., \cite[Lemma 3.6]{CP}) that \[\wh{U}\cong\mc{M}(G,v,v,I)\] with $G$ a maximal subgroup of $\wh{U}$ (that is this isomorphic to
$\gr(S\cap G_1)$) and $G$ is abelian-by-finite (we denote by $I$ the identity matrix of degree $v$). Put $\wh{S}=(S\setminus U)\cup \wh{U}$,
a disjoint union. Note that $\wh{S}$ also is a subsemigroup of $\M_m(L)$ and $\wh{U}$ is an ideal of $\wh{S}$ (cf. \cite[Lemma 2.5.2]{JO}).
Hence, $K_0[S]$ is a subalgebra of $K_0[\wh{S}]$ and it has $K_0[\wh{U}]$ as an ideal. The ideal $K_0[\wh{U}]$ has $e=e_1+\dotsb+e_v$
as an identity and thus this is a central element of $K_0[\wh{S}]$. We also have a natural epimorphism \[f_P\colon K_0[S]\to K_0[S]e\s K_0[\wh{U}].\]
Hence, $K_0[S]e$ is a Noetherian algebra and $K_0[S]e\s K_0[\wh{U}]\cong\M_v(K[G])$. By \cite[Proposition 2.5.6]{JO} we also know that $G$
is finitely generated. So, $G$ is a finitely generated abelian-by-finite group. Note that
\[\Ker f_P=\{\alpha\in K_0[S]:\alpha e=0\}=\{\alpha\in K_0[S]:\alpha U=0\}.\]
Since the ideal $K_0[U]$ is not contained in the prime ideal $P/K[P\cap M]$, we get $\Ker f_P\s P/K[P\cap M]$.

We are now in a position to prove the main result of this section
(in the statement and proof of this result we use the notation introduced in this section).

\begin{thm}
	If $P$ is a prime ideal of $K[M]$ then there exists an ideal $I_P$ of $K[M]$ contained in $P$ and a prime ideal $P_0$ of $K[G]$ such that
	$K[M]/I_P\s\M_v(K[G])$ and $K[M]/P\s\M_v(K[G]/P_0)$ for some $v\ge 1$. Moreover, $\M_v(K[G])$ is a localization of $K[M]/I_P$.
	In particular, $\Qcl(K[M]/P)\cong\Qcl(\M_v(K[G]/P_0))$. If, furthermore, $K[M]$ is semiprime then there exist finitely
	many finitely generated abelian-by-finite groups, say $G_1,\dotsc,G_m$, each being the group of quotients of a cancellative
	subsemigroup of $M$, such that $K[M]$ embeds into $\M_{v_1}(K[G_1])\times\dotsb\times\M_{v_m}(K[G_m])$ for some $v_1,\dotsc,v_m\ge 1$.
	\begin{proof}
		Let $P$ be a prime ideal of $K[M]$. If $P\cap M=\vn$ then the first part of the result has been shown in Proposition \ref{prop:11}.
		So, assume that $P\cap M\ne\vn$. Let $S=M/(P\cap M)$. From the above we know that $K_0[S]/\Ker f_P\s K_0[\wh{U}]\cong\M_v(K[G])$,
		where $G$ is the group of fractions of a cancellative subsemigroup of $U$.
		Furthermore, $K_0[\wh{U}] $ is a localization of $K_0[U]$ with respect to diagonal matrices (with entries in $G$) that belong
		to $K_0[U]$. Such matrices are regular in $K_0[\wh{U}]$ and thus they also are regular elements in $K_0[S]/\Ker f_P$.
		Hence the Noetherian algebra $K_0[\wh{U}]$ is a localization of $K_0[S]/\Ker f_P$. Therefore, as it is well-known
		(see \cite[Theorem 3.2.6]{JO}), there exists a prime ideal $P'_0=\M_v(P_0)$ of $\M_v(K[G])$ (with $P_0$
		a prime  ideal of $K[G]$) such that $P'_0\cap(K_0[S]/\Ker f_P)=P/\Ker f_P$. Let $I_P$ denote the ideal of $K[M]$ containing
		$K[P\cap M]$ that naturally projects onto $\Ker f_P$ in $K_0[S]$. It follows that $K[M]/I_P \s\M_v(K[G])$ and
		$K[M]/P\s\M_v(K[G]/P_0)$ and $K[M]/P$ is a localization of $\M_v(K[G]/P_0)$. Hence the first part of the result follows.
		
		Assume now the the algebra $K[M]$ is semiprime. Because $K[M]$ is Noetherian (see Theorem \ref{thm:5}), it has finitely many minimal
		prime ideals, say $P_1,\dotsc,P_m$. By the first part, for each $P_i$ there exists an ideal $I_{P_i}\s P_i$ such that 
		$K[M]/I_{P_i}\s\M_{v_i}(K[G_i])$, for some finitely generated abelian-by-finite group $G_i$ that is the group of fractions of
		a cancellative subsemigroup of $M$. Since $\bigcap_{i=1}^mI_i\s\bigcap_{i=1}^mP_i=0$, we get that $K[M]$ embeds into 
		$K[M]/I_{P_1}\times\dotsb\times K[M]/I_{P_m}$. Hence the result follows.
	\end{proof}
\end{thm}

\bibliographystyle{amsplain}
\bibliography{refs}

\end{document}